\documentclass[twoside]{bourlarge}
\usepackage{modele_h}
\usepackage{amssymb}
\usepackage{amsmath}
\usepackage[english]{babel}
\usepackage{graphicx}

\newcommand{\nwc}{\newcommand}
\nwc{\nwt}{\newtheorem}
\nwt{coro}{Corollary}
\nwt{ex}{Example}
\nwt{rem}{Remark}
\nwt{Theorem}{Theorem}
\nwt{Conj}{Conjecture}
\nwt{Question}{Question}

\nwc{\mf}{\mathbf} 
\nwc{\blds}{\boldsymbol} 
\nwc{\ml}{\mathcal} 

\nwc{\lam}{\lambda}
\nwc{\del}{\delta}
\nwc{\Del}{\Delta}
\nwc{\Lam}{\Lambda}
\nwc{\elll}{\ell}

\nwc{\IA}{\mathbb{A}} 
\nwc{\IB}{\mathbb{B}} 
\nwc{\IC}{\mathbb{C}} 
\nwc{\ID}{\mathbb{D}} 
\nwc{\IE}{\mathbb{E}} 
\nwc{\IF}{\mathbb{F}} 
\nwc{\IG}{\mathbb{G}} 
\nwc{\IH}{\mathbb{H}} 
\nwc{\IN}{\mathbb{N}} 
\nwc{\IP}{\mathbb{P}} 
\nwc{\IQ}{\mathbb{Q}} 
\nwc{\IR}{\mathbb{R}} 
\nwc{\IS}{\mathbb{S}} 
\nwc{\IT}{\mathbb{T}} 
\nwc{\IZ}{\mathbb{Z}} 

\def\bbleft{{\mathchoice {[\mskip-3mu {[}} {[\mskip-3mu {[}}{[\mskip-4mu {[}}{[\mskip-5mu {[}}}}
\def\bbright{{\mathchoice {]\mskip-3mu {]}} {]\mskip-3mu {]}}{]\mskip-4mu {]}}{]\mskip-5mu {]}}}}
\nwc{\setK}{\bbleft 1,K \bbright}
\nwc{\setN}{\bbleft 1,\cN \bbright}

\nwc{\va}{{\bf a}}
\nwc{\vb}{{\bf b}}
\nwc{\vc}{{\bf c}}
\nwc{\vd}{{\bf d}}
\nwc{\ve}{{\bf e}}
\nwc{\vf}{{\bf f}}
\nwc{\vg}{{\bf g}}
\nwc{\vh}{{\bf h}}
\nwc{\vi}{{\bf i}}
\nwc{\vI}{{\bf I}}
\nwc{\vj}{{\bf j}}
\nwc{\vk}{{\bf k}}
\nwc{\vl}{{\bf l}}
\nwc{\vm}{{\bf m}}
\nwc{\vM}{{\bf M}}
\nwc{\vn}{{\bf n}}
\nwc{\vo}{{\it o}}
\nwc{\vp}{{\bf p}}
\nwc{\vq}{{\bf q}}
\nwc{\vr}{{\bf r}}
\nwc{\vs}{{\bf s}}
\nwc{\vt}{{\bf t}}
\nwc{\vu}{{\bf u}}
\nwc{\vv}{{\bf v}}
\nwc{\vw}{{\bf w}}
\nwc{\vx}{{\bf x}}
\nwc{\vy}{{\bf y}}
\nwc{\vz}{{\bf z}}
\nwc{\bal}{\blds{\alpha}}
\nwc{\bep}{\blds{\epsilon}}
\nwc{\barbep}{\overline{\blds{\epsilon}}}
\nwc{\bnu}{\blds{\nu}}
\nwc{\bmu}{\blds{\mu}}

\nwc{\bk}{\blds{k}}
\nwc{\bm}{\blds{m}}
\nwc{\bM}{\blds{M}}
\nwc{\bp}{\blds{p}}
\nwc{\bq}{\blds{q}}
\nwc{\bn}{\blds{n}}
\nwc{\bv}{\blds{v}}
\nwc{\bw}{\blds{w}}
\nwc{\bx}{\blds{x}}
\nwc{\bxi}{\blds{\xi}}
\nwc{\by}{\blds{y}}
\nwc{\bz}{\blds{z}}

\nwc{\cA}{\ml{A}}
\nwc{\cB}{\ml{B}}
\nwc{\cC}{\ml{C}}
\nwc{\cD}{\ml{D}}
\nwc{\cE}{\ml{E}}
\nwc{\cF}{\ml{F}}
\nwc{\cG}{\ml{G}}
\nwc{\cH}{\ml{H}}
\nwc{\cI}{\ml{I}}
\nwc{\cJ}{\ml{J}}
\nwc{\cK}{\ml{K}}
\nwc{\cL}{\ml{L}}
\nwc{\cM}{\ml{M}}
\nwc{\cN}{\ml{N}}
\nwc{\cO}{\ml{O}}
\nwc{\cP}{\ml{P}}
\nwc{\cQ}{\ml{Q}}
\nwc{\cR}{\ml{R}}
\nwc{\cS}{\ml{S}}
\nwc{\cT}{\ml{T}}
\nwc{\cU}{\ml{U}}
\nwc{\cV}{\ml{V}}
\nwc{\cW}{\ml{W}}
\nwc{\cX}{\ml{X}}
\nwc{\cY}{\ml{Y}}
\nwc{\cZ}{\ml{Z}}

\nwc{\tA}{\widetilde{A}}
\nwc{\tB}{\widetilde{B}}
\nwc{\tE}{E^{\vareps}}
\nwc{\tk}{\tilde k}
\nwc{\tN}{\tilde N}
\nwc{\tP}{\widetilde{P}}
\nwc{\tQ}{\widetilde{Q}}
\nwc{\tR}{\widetilde{R}}
\nwc{\tV}{\widetilde{V}}
\nwc{\tW}{\widetilde{W}}
\nwc{\ty}{\tilde y}
\nwc{\teta}{\tilde \eta}
\nwc{\tdelta}{\tilde \delta}
\nwc{\tlambda}{\tilde \lambda}
\nwc{\ttheta}{\tilde \theta}
\nwc{\tvartheta}{\tilde \vartheta}
\nwc{\tPhi}{\widetilde \Phi}
\nwc{\tpsi}{\tilde \psi}

\nwc{\To}{\longrightarrow} 

\nwc{\ad}{\rm ad}
\nwc{\eps}{\epsilon}
\nwc{\ep}{\epsilon}
\nwc{\vareps}{\varepsilon}

\nwc{\err}{\epsilon}
\def\ep{\epsilon}

\def\sq2{\sqrt{2}}

\def\t2{{\mathbb T}^2}
\def\s2{{\mathbb S}^2}
\def\hn{\mathcal{H}_{N}}

\nwc{\lap}{\bigtriangleup}
\nwc{\rest}{\restriction}
\nwc{\Diff}{\operatorname{Diff}}
\nwc{\diam}{\operatorname{diam}}
\nwc{\Res}{\operatorname{Res}}
\nwc{\Spec}{\operatorname{Spec}}
\nwc{\Vol}{\operatorname{Vol}}
\nwc{\Var}{\operatorname{Var}}
\nwc{\Op}{\operatorname{Op}}
\nwc{\Oph}{\operatorname{Op}_\hbar}
\nwc{\supp}{\operatorname{supp}}
\nwc{\Span}{\operatorname{span}}

\nwc{\dia}{\varepsilon}
\nwc{\cut}{f}
\nwc{\qm}{u_\hbar}

\def\hto0{\xrightarrow{\hbar\to 0}}

\def\rto0{\xrightarrow{r\to 0}}

\nwc{\la}{\langle}
\nwc{\ra}{\rangle}
\nwc{\lp}{\left(}
\nwc{\rp}{\right)}

\nwc{\bequ}{\begin{equation}}
\nwc{\be}{\begin{equation}}
\nwc{\ben}{\begin{equation*}}
\nwc{\bea}{\begin{eqnarray}}
\nwc{\bean}{\begin{eqnarray*}}
\nwc{\bit}{\begin{itemize}}
\nwc{\bver}{\begin{verbatim}}

\nwc{\eequ}{\end{equation}}
\nwc{\ee}{\end{equation}}
\nwc{\een}{\end{equation*}}
\nwc{\eea}{\end{eqnarray}}
\nwc{\eean}{\end{eqnarray*}}
\nwc{\eit}{\end{itemize}}
\nwc{\ever}{\end{verbatim}}

\newcommand{\defeq}{\stackrel{\rm{def}}{=}}

\begin{document}
\Large
\title{Anatomy of quantum chaotic eigenstates}

\author{St\'ephane {\sc Nonnenmacher} \footnote{I am grateful to
    E.Bogomolny, who allowed me to reproduce several plots from
    \cite{BogoSchm07}. The author has been partially supported
by the Agence Nationale de la Recherche under the grant ANR-09-JCJC-0099-01.
These notes were written while he was visiting the Institute of Advanced Study in Princeton,
supported by the National Science Foundation under agreement No. DMS-0635607.}\\ 
Institut de Physique Th\'eorique\\ 
CEA-Saclay\\
91191 Gif-sur-Yvette, France}


\maketitle

\section{Introduction}\label{sec1}

These notes present a description of {\em quantum chaotic
eigenstates}, that is bound states of a quantum dynamical system, the classical
limit of which is chaotic. The classical dynamical systems we will be dealing
with are mostly of two
types: geodesic flows on euclidean domains (``billiards'') or compact
riemannian manifolds, and canonical
transformations on a compact phase space; the common feature is the
``chaoticity'' of the dynamics. The corresponding quantum systems will
always be considered in the semiclassical (or
high-frequency) r\'egime, in order to establish a link with the classical
dynamics. As a first illustration, we plot below two eigenstates of a paradigmatic system,
the Laplacian on the {\em stadium billiard}, with Dirichlet boundary conditions\footnote{The
  eigenfunctions of the stadium plotted in this article were computed using a code
  nicely provided to me by E.~Vergini, which uses the {\em scaling
    method} invented in \cite{VS95}.}. 

The study of chaotic eigenstates makes up a large part of the field of
{\em quantum chaos}. It is somewhat
complementary with the contribution of J.Keating, who will focus on
the statistical properties of quantum spectra, another major topic
in quantum chaos. I do not include the study of eigenstates of
quantum graphs (a recent interesting development in the field), since
this question should be addressed in
U.Smilansky's lecture. Although these notes are purely theoretical,
H.-J.~St\"ockmann's lecture will show that the questions raised here have
direct experimental applications (his lecture should present
experimentally measured eigenmodes of 2- and 3-dimensional ``quantum billiards'').
\begin{figure}[htbp]
\begin{center}
\includegraphics[angle=-90,width=1.\textwidth]{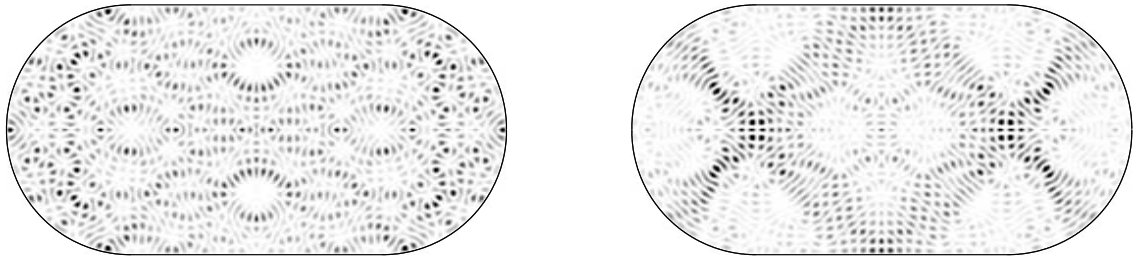}\hspace{1cm}
\end{center}
\caption{Two eigenfunctions of the Dirichlet Laplacian in the stadium
  billiard, with wavevectors $k=60.196$ and $k=60.220$ (see
  \eqref{e:Helmholtz}). Large values of $|\psi(x)|^2$ correspond to
  dark regions, while nodal lines are white. While the left eigenfunction
  looks relatively ``equidistributed'', the right one is {\em scarred} by two
  symmetric periodic orbits (see \S\ref{s:scars}).\label{f:stadium-eigen}}
\end{figure}

One common feature of the chaotic eigenfunctions (except in some very
specific systems) is the absence of explicit
(or even approximate) formulas for them. One then has to resort to
indirect, rather unprecise approaches to describe these
eigenstates. We will use various analytic tools or points of view: 
deterministic/statistical, macro/microscopic, pointwise/global properties,
generic/specific systems. The
level of rigour in the results varies from mathematical
proofs to heuristics, generally supported by numerical
experiments. 
The necessary selection of results reflects my
personal view or knowledge of the subject, it omits
several important developments, and is more ``historical'' than
sharply up-to-date. The list of references is qisw but in no way exhaustive. 

These notes are organized as follows. We introduce in section~\ref{s:chaos}
the classical dynamical systems we will focus on (mostly geodesic
flows and maps on the 2-dimensional torus), mentioning their
degree of ``chaos''. We also sketch the 
quantization procedures leading to quantum Hamiltonians or
propagators, the eigenstates of which we want to understand. We also
mention some properties of the semiclassical/high-frequency limit.
In section~\ref{s:macro} we describe the {\em macroscopic} properties of
the eigenstates in the semiclassical limit, embodied by their
{\em semiclassical measures}. These
properties include the {\em quantum ergodicity} property,
which for some systems (with arithmetic symmetries) can be improved to
{\em quantum unique ergodicity}, namely the fact that all high-frequency eigenstates
are ``flat'' at the macroscopic scale; on the opposite, some specific systems allow the presence of
exceptionally localized eigenstates. In section~\ref{s:statistics} we focus
on more refined properties of the eigenstates, many of {\em statistical} nature (value
distribution, correlation functions). Very little is
known rigorously, so one has to resort to heuristical models of {\em random wavefunctions} to
describe these statistical properties. The large values of the
wavefunctions or Husimi densities are discussed, including the {\em
  scar phenomenon}. Section~\ref{s:nodal} discusses the most ``quantum''
or microscopic aspect of the eigenstates, namely their {\em nodal sets},
both in position and phase space (Husimi) representations. Here as
well, the random state models are helpful, and lead to interesting
questions in probability theory.

\section{What is a quantum chaotic eigenstate?\label{s:chaos}}

In this section we first present a general definition of the notion of
``chaotic eigenstate''. We then focus our attention to geodesic flows
in euclidean domains or on compact riemannian manifolds of negative curvature, which form
the simplest systems proved to be chaotic. Finally we present some
discrete time dynamics (chaotic canonical maps on the 2-dimensional torus).

\subsection{A short review of quantum mechanics} 

Let us start by recalling that classical mechanics on the phase space $T^*\IR^d$
can be defined, in the Hamiltonian formalism, by a real valued function $H(x,p)$
on that phase space, called the Hamiltonian. We will always assume the
system to be autonomous, namely the function $H$ to be independent of
time. This function then generates
the flow\footnote{We always assume that the flow is complete, that is it does
not blow up in finite time.}
$$
\vx(t)\defeq (x(t),p(t))=\Phi^t_H(\vx(0)),\qquad t\in\IR\,,\quad
\vx(t)\in T^*\IR^d\,,
$$ 
by solving Hamilton's equations (as usual $\dot{f}$ is the time derivative):
\be\label{eq:Hamilton}
\dot x_j(t)= \frac{\partial H}{\partial p_j}(\vx(t)),\quad 
\dot p_j(t)= -\frac{\partial H}{\partial x_j}(\vx(t)),\quad j=1,\ldots,d\,.
\ee
This flow preserves the symplectic form $\sum_jdp_j\wedge dx_j$, and
the energy shells $\cE_E=H^{-1}(E)$.

The corresponding quantum
mechanical system is defined by an operator $\hat H_\hbar$ acting on
the (quantum) Hilbert space $\cH=L^2(\IR^d,dx)$. This operator can be
formally obtained by replacing coordinates $x,p$ by operators:
\be\label{e:hatH}
\hat H_\hbar=H (\hat x_\hbar,\hat p_\hbar),
\ee 
where $\hat x_\hbar$ is the operator of multiplication by $x$, while the
momentum operator $\hat p_\hbar=\frac{\hbar}{i}\nabla$ 
is conjugate to $\hat x$ through the $\hbar$-Fourier
transform $\cF_\hbar$.
The notation \eqref{e:hatH} assumes that one has selected a certain ordering
between the operators $\hat x_\hbar$ and $\hat p_\hbar$;
in physics one usually chooses the fully symmetric ordering, also
called the {\em Weyl
quantization}: it has the advantage to make $\hat H_\hbar$ formally
self-adjoint on $L^2(\IR^d)$. Quantization procedures can
also be defined when the euclidean space $\IR^d$ is
replaced by a compact manifold $M$, by gluing together quantizations
defined on different coordinate patches.

The quantum dynamics, which governs the evolution of the wavefunction
$\psi(t)\in\cH$ describing the system, 
is then given by the Schr\"odinger equation:
\be\label{e:Schroe}
i \hbar\, \frac{\partial \psi(t)}{\partial t}=[\hat H_\hbar\, \psi](t)\,.
\ee
Solving this linear equation produces the propagator, that is the
family of unitary operators on $L^2(\IR^d)$,
$$
U_\hbar^t = \exp(-i\hat H_\hbar/\hbar)\,,\quad t\in\IR\,.
$$
\begin{rem}In physical systems, Planck's constant $\hbar$ is a fixed number,
which is of order $10^{-34}$ in standard units. However, if the
system (atom, molecule, ``quantum dot'') is itself microscopic, the
value of $\hbar$ may be comparable with the {\em typical action} of the
system, in which case it is more natural to select units in which
$\hbar=1$. 
Our point of view throughout this work will be the opposite: we will
assume that $\hbar$ is (very) small
compared with the typical action of the system, and most results will
be valid asymptotically, in the {\em semiclassical limit} $\hbar\to 0$. 
\end{rem}

\subsection{Quantum-classical correspondence}
At this point, let us 
introduce the crucial semiclassical property of the quantum
evolution: it is called (in the physics literature) the
quantum-classical correspondence, while in mathematics this result is
known as Egorov's theorem. This property states that the evolution of
{\em observables}
approximately commutes with quantization.
For us, a (classical) observable is a smooth, compactly supported function on phase space $f\in
C_c^\infty(T^*\IR^d)$, while the corresponding quantum observable
$\hat{f}_\hbar$ is obtained through the above mentioned quantization
procedure. The evolution of classical and quantum
evolutions are defined by duality with that of
particles/wavefunctions:
$$
f(t)=f\circ \Phi^t_H,\qquad \hat f_\hbar (t)=U_\hbar^{-t}\hat f_\hbar U_\hbar^t\,.
$$ 
The quantum-classical
correspondence connects these two evolutions:
\be\label{e:Egorov}
\forall t\in\IR,\qquad \hat f_\hbar (t)=
\widehat{f(t)}_\hbar
+\cO_{L^2\to L^2}(e^{\Gamma|t|}\hbar)\,,\qquad
\ee
where the exponent $\Gamma>0$ depends on the instability of the
classical flow.

\medskip

The most common form of dynamics is the motion of a scalar
particle in an electric potential $V(x)$. It corresponds to the Hamiltonian
\be\label{eq:Schrod}
H(x,p)=\frac{|p|^2}{2m} + qV(x),\quad\text{quantized into}\quad \hat H_\hbar =
-\frac{\hbar^2\Delta}{2m} + qV(x)\,.
\ee
We will usually scale the mass and electric charge to $m=q=1$, keeping
$\hbar$ small. 
Since the Hamilton flow \eqref{eq:Hamilton} leaves each energy energy
shell $\cE_E$ invariant, we may restrict our attention to the flow on
a single shell. We will be interested in cases where 
\begin{enumerate}
\item the energy shell $\cE_E$ is {\em bounded} in phase space (that is, both
  the positions and momenta of the particles remain finite at all
  times). This is the case if $V(x)$ is confining ($V(x)\to\infty$ as $|x|\to\infty$).
\item the flow on $\cE_E$ is {\em chaotic}
  (and this is also the case on the neighbouring shells $\cE_{E+\eps}$). 
  ``Chaos'' is a vague word, which we will make more precise below. 
\end{enumerate}
The first condition implies that, provided $\hbar$ is small enough, the
spectrum of $\hat H_\hbar$ is purely discrete near the energy $E$,
with eigenstates $\psi_{\hbar,j}\in L^2(\IR^d)$ ({\em
  bound states}). Besides, fixing some small $\eps>0$ and
letting $\hbar\to 0$, the number of 
eigenvalues $E_{\hbar,j}\in [E-\eps,E+\eps]$ typically grows like
$C \hbar^{-d}$.
Under the second condition, the eigenstates
with energies in this interval can be called {\em quantum chaotic eigenstates}. 

Below we describe several degrees of ``chaos'', which regard the {\em long time properties} of the classical
flow. These properties are relevant when
describing the eigenstates of the quantum system, which
form the ``backbone'' of the long time quantum dynamics. 
The main objective of quantum chaos consists in connecting, in
a precise way, the classical and quantum long time (or time
independent) properties.

\subsection{Various levels of chaos}
For most Hamiltonians of the form \eqref{eq:Schrod}
(e.g. the physically relevant case of a hydrogen atom in a constant
magnetic field),
the classical dynamics on bounded energy shells $\cE_E$
involves both regular and chaotic regions of phase space; one
then speaks of a {\em mixed dynamics} on $\cE_E$. The regular region is composed of a
number of ``islands of stability'', made of quasiperiodic motion structured around
stable periodic orbits; these islands are embedded in a
``chaotic sea'' where trajectories are unstable (they have a
positive Lyapunov exponent). These notions
of ``island of stability'' versus ``chaotic sea'' are rather
poorly understood mathematically, but have received compelling
numerical evidence \cite{LichLieb92}. The
main conjecture concerning the corresponding quantum system, is that
most eigenstates are either localized in the regular region, or in the chaotic
sea \cite{Percival}. To my knowledge this conjecture remains fully open at present,
in part due to our lack of understanding of the classical dynamics. 

For this reason, I will restrict myself (as most researches in quantum
chaos do) to the case of systems admitting a {\em purely chaotic} dynamics
on $\cE_E$. I will allow various degrees of chaos, the minimal
assumption being the {\em ergodicity} of the flow $\Phi^t_H$ on $\cE_E$, with
respect to the natural (Liouville) measure on $\cE_E$. This assumption
means that, for almost every initial point $\vx_0\in \cE_E$, the
time averages of any observable $f$ converge to its phase
space average:
\be\label{e:ergodicity}
\lim_{T\to\infty}\frac{1}{2T}\int_{-T}^T f(\Phi^t_H(\vx_0))\,dt =
\int_{\cE_E}f(\vx)\,d\mu_L(\vx) \defeq \frac{\int f(\vx)\,\delta(H(\vx)-E)\,d\vx}{\int \delta(H(\vx)-E)\,d\vx}\,.
\ee
A stronger chaotic property is the {\em mixing} property,
or decay of time correlations between two observables $f,g$:
\be\label{e:mixing}
C_{f,g}(t)\defeq \int_{\cE_{E}} g\times (f\circ\Phi^t)\,d\mu_L -
\int f\,d\mu_L\int g\,d\mu_L \xrightarrow{t\to\infty}0\,.
\ee
The decay rate of $C_{f,g}(t)$ (or {\em mixing rate})
depends on both the flow $\Phi^t$ and the regularity of the
observables $f,g$. For very chaotic flows (Anosov flows, see \S\ref{s:geodesic}) and
smooth observables, the decay is exponential.

\subsection{Geometric quantum chaos}

In this section we give explicit examples of chaotic flows, namely the
geodesic flows in certain euclidean billiards compact riemannian manifolds. 
The dynamics is then induced by the geometry, rather than by a
potential. Both the classical and quantum properties of these systems
have been investigated a lot in the past 30 years.

\subsubsection{Billiards}
\begin{figure}[htbp]
\begin{center}
\includegraphics[width=1.\textwidth]{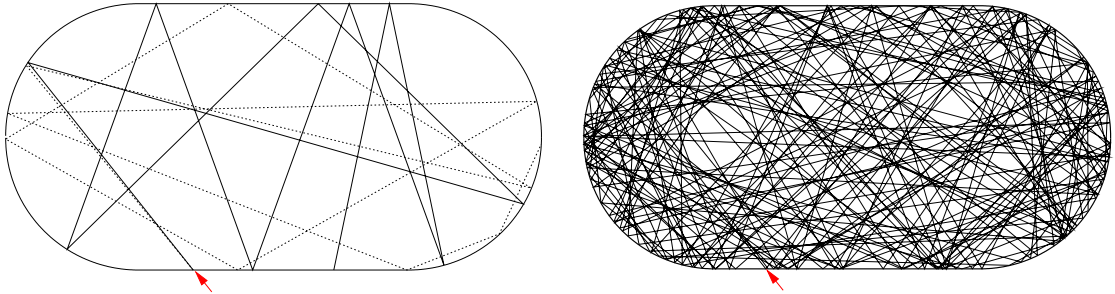}
\end{center}
\caption{Left: two trajectories in the stadium billiard, initially
  very close to one another, and then diverging fast (the
  red arrow shows the initial point). Right:
long evolution of one of the trajectories.}\label{f:stadium}
\end{figure}

The simplest form of ergodic system occurs when the potential $V(x)$
is an infinite barrier delimiting a bounded domain
$\Omega\subset\IR^d$ (say, with piecewise smooth boundary), so that the particle moves freely inside
$\Omega$ and bounces specularly at the boundaries. For obvious
reasons, such a system is called a {\em euclidean billiard}. All positive energy
shells are equivalent to one another, up to a
rescaling of the velocity, so we may restrict our attention to the
shell $\cE=\{(x,p),\,x\in\Omega,\,|p|=1\}$. The long
time dynamical properties only depend on the shape of the
domain. 
For instance, in 2 dimensions, rectangular, circular or
elliptic billiards lead to
{\em integrable} dynamics: the flow admits
two independent integrals of motion --- in the case of
the circle, the energy and the angular momentum. A
convex billiard with a smooth boundary will always admit
``whispering gallery'' stable orbits, so such billiards cannot be
fully chaotic. On the opposite, the famous {\em
  stadium billiard} (see Fig.~\ref{f:stadium})
  was proved to be ergodic by Bunimovich \cite{Bunim}.
Historically, the first euclidean billiard proved to be ergodic was the Sinai billiard,
composed of one or several circular obstacles
inside a square (or torus) \cite{Sinai70}. 
These billiards
also have positive Lyapunov exponents (meaning that almost
all trajectories are exponentially unstable, see the left part of Fig.~\ref{f:stadium}). 
It has been shown more recently that these billiards
are {\em mixing}, but with correlations decaying at polynomial or
subexponential rates \cite{Chernov07,BalMel08,Melb09}. 

The quantization of the broken geodesic flow inside $\Omega$ is the
(semiclassical) Laplacian with Dirichlet boundary conditions:
\be\label{e:Laplacian}
\hat H_\hbar = -\frac{\hbar^2\Delta_\Omega}{2}\,.
\ee
Obviously, the parameter $\hbar$ only amounts to a rescaling of the
spectrum:
an eigenstate $\psi_\hbar$ of \eqref{e:Laplacian} with energy
$E_\hbar\approx 1/2$ is also an eigenstate
of $-\Delta_\Omega$ with eigenvalue $k^2\approx \hbar^{-2}$. Hence, $\hbar$
represents the {\em wavelength} of $\psi_\hbar$, the inverse of its
{\em wavevector} $k$. 
Fixing $E=1/2$ and taking the semiclassical limit $\hbar\to 0$ is equivalent
with studying the high-frequency or high-wavevector spectrum of $-\Delta_\Omega$.

The system \eqref{e:Laplacian} is often called a {\em quantum
  billiard}, although this operator is not only relevant in quantum
mechanics, but in all sorts of wave mechanics (see H.-J.~St\"ockmann's
lecture). Indeed, the scalar Helmholtz equation
\be\label{e:Helmholtz}
\Delta \psi_j+k_j^2\psi_j=0\,,
\ee 
may describe stationary acoustic waves in a cavity, and is also relevant to
describe electromagnetic waves in a quasi-2D cavity, provided one is
allowed to separate the different polarization components of
the electric field. 

Euclidean billiards thus form the simplest {\em
  realistic} quantized chaotic systems, for which the classical
dynamics is well understood at the mathematical level. Besides, 
the spectrum of the Dirichlet Laplacian can be numerically computed
up to large values of $k$ using methods specific to
the euclidean geometry, like the scaling method \cite{VS95}. 
For these reasons, these billiards have become 
a paradigm of quantum chaos studies.

\subsubsection{Anosov geodesic flows\label{s:geodesic}}
\begin{figure}[htbp]
\begin{center}
\includegraphics[width=.5\textwidth]{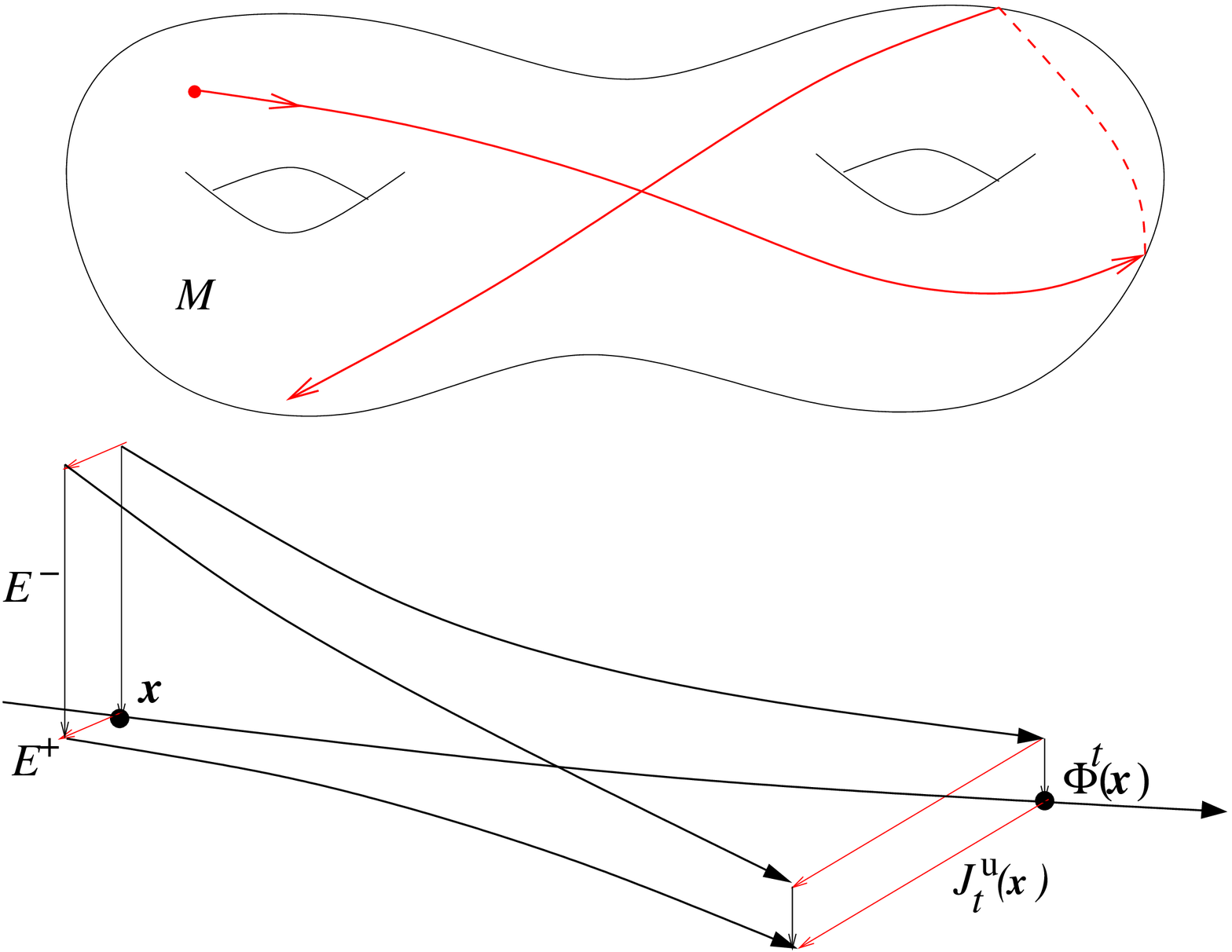}\hspace{1cm}
\includegraphics[width=.35\textwidth]{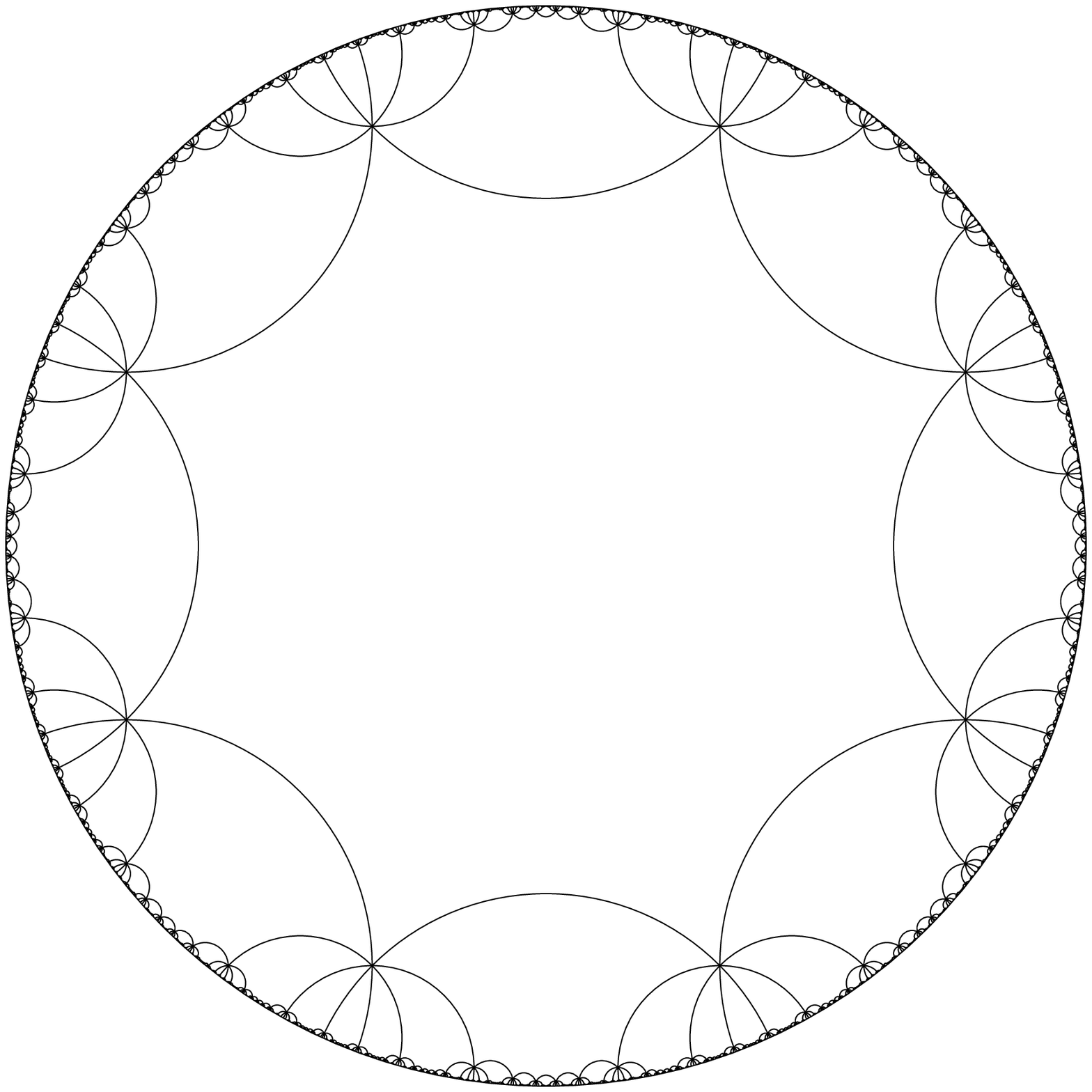}
\end{center}
\caption{Top left: geodesic flow on a surface of negative
  curvature (such a surface has a genus $\geq 2$). Right:
  fundamental domain for an ``octagon'' surface $\Gamma\backslash\IH$
  of constant negative curvature (the figure is due to C.~McMullen)
Botton left: a phase space trajectory and two nearby trajectories
approaching it in the future or past. The stable/unstable directions at
$\vx$ and $\Phi^t(\vx)$ are shown. The red lines feature the expansion along the
unstable direction, measured by the
unstable Jacobian $J_t^u(\vx)=\det(d\Phi^t\rest_{E^+_\vx})$.}\label{f:manifolds}
\end{figure}

The strongest form of chaos occus in
systems (maps or flows) with the {\em Anosov property}, also called
{\em uniformly hyperbolic systems} \cite{Anosov67}. The first (and
main) example of an Anosov flow is given by the {\em geodesic flow on
  a compact riemannian manifold $(M,g)$ of negative curvature},
generated by the free particle Hamiltonian $H(x,p)=|p|^2_g/2$.
Uniform hyperbolicity --- which is induced by the negative
curvature of the manifold --- means that at each point
$\vx\in\cE$ the tangent space $T_\vx\cE$ splits into the vector
$X_\vx$ generating the flow, the unstable subspace $E^+_\vx$ and the stable subspace $E^-_\vx$.
The stable (resp. unstable) subspace is defined by the property that the
flow contracts vectors exponentially in the future (resp. in the past),
see Fig.~\ref{f:manifolds}: for some $C>0,\ \lambda>0$ independent of $\bx\in\cE$,
\be\label{e:Anosov}
\forall v\in E^\pm_\vx,\ \forall t>0,\quad \|d\Phi^{\mp t}_\vx v\|\leq C\,e^{-\lambda t}\|v\|.
\ee
Anosov systems present the strongest form of chaos, but their ergodic
properties are (paradoxically) better understood than for
the billiards of the previous section. The flow has a positive {\em
  complexity}, reflected in the exponential proliferation of long periodic geodesics.
For this reason, this geometric model has been at the center of the mathematical
investigations of quantum chaos, in spite of its minor physical relevance. 
Generalizing the case of the euclidean billiards, the quantization of the
geodesic flow on $(M,g)$ is given by the (semiclassical)
Laplace-Beltrami operator
\be
\hat H_\hbar = -\frac{\hbar^2\Delta_g}{2}
\ee
acting on the Hilbert space $L^2(M,dg)$ associated with the
Lebesgue measure. The eigenstates of $\hat H_\hbar$ with eigenvalues
$\approx 1/2$ (equivalently, the high-frequency eigenstates of
$-\Delta_g$) constitute a class of quantum chaotic
eigenstates, the study of which is not impeded by boundary
problems present in billiards.

The spectral properties of this Laplacian have interested
mathematicians working in riemannian geometry, PDEs, analytic number
theory, representation theory, for at least a century, while the specific
``quantum chaotic'' aspects have emerged only in the last 30 years.

The first example of a manifold with negative curvature
is the Poincar\'e half-space $\IH$  with its hyperbolic
metric $\frac{dx^2+dy^2}{y^2}$, on which
the group $SL_2(\IR)$ acts isometrically by Moebius
transformations  ($\IH$ is often represented by the
equivalent {\em hyperbolic disk}, like in Fig.~\ref{f:manifolds}).
For certain discrete subgroups $\Gamma$ of  $PSL_2(\IR)$
(called co-compact lattices), the quotient $M=\Gamma\backslash \IH$ is a
smooth compact surface. This group structure provides detailed
information on the spectrum of the Laplacian (for instance, the
Selberg trace formula explicitly connects the spectrum of $\Delta_{g}$
with the periodic geodesics on $M$, which are themselves represented by the conjugacy classes
of $\Gamma$). 

Furthermore, for some of these discrete subgroups $\Gamma$ (called
{\em arithmetic}), one can construct a {\em commutative} algebra of
Hecke operators on $L^2(M)$, which also commute with the Laplacian; it
then make sense to study in priority the joint eigenstates of $\Delta$
and of these Hecke operators, which we will call the Hecke eigenstates.
This arithmetic structure allows to obtain extra nontrivial information on these
eigenstates (see \S\ref{s:arithmetic}), so the latter will
appear several times along these notes. Their study composes a part of {\em arithmetic quantum
  chaos}, a lively field of research.

\subsection{Classical and quantum chaotic maps}
Beside the Hamiltonian or geodesic flows, another model system has
attracted much attention in the dynamical systems community: chaotic
maps on some compact phase space $\cP$. Instead of a flow, the
dynamics is given by a discrete
time transformation $\kappa:\cP\to\cP$. Because we want to quantize
these maps, we require the phase space $\cP$ to have a symplectic
structure, and the map $\kappa$ to preserve this structure (in other words,
$\kappa$ is an invertible canonical transformation on $\cP$). 
\begin{figure}[htbp]
\begin{center}
\includegraphics[width=.6\textwidth]{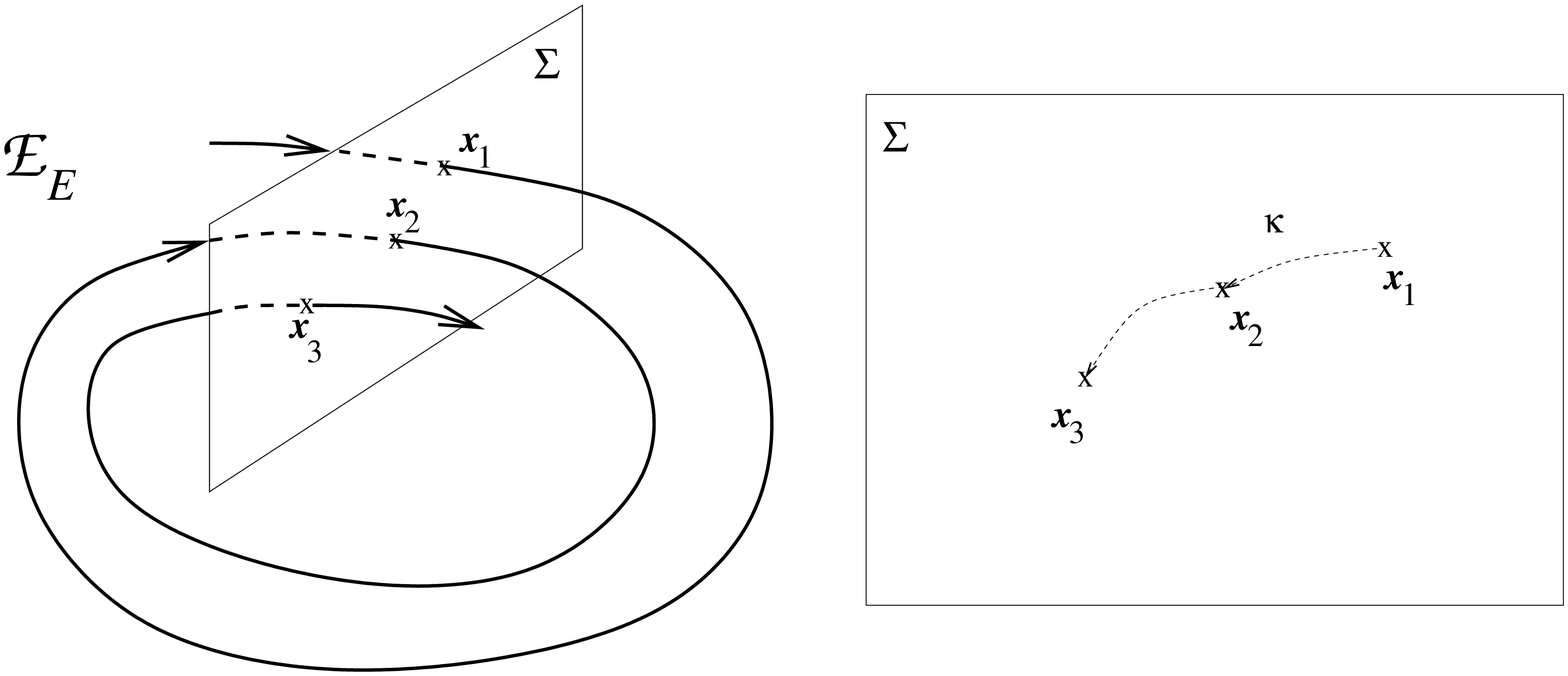}\hspace{1cm}
\includegraphics[width=.6\textwidth]{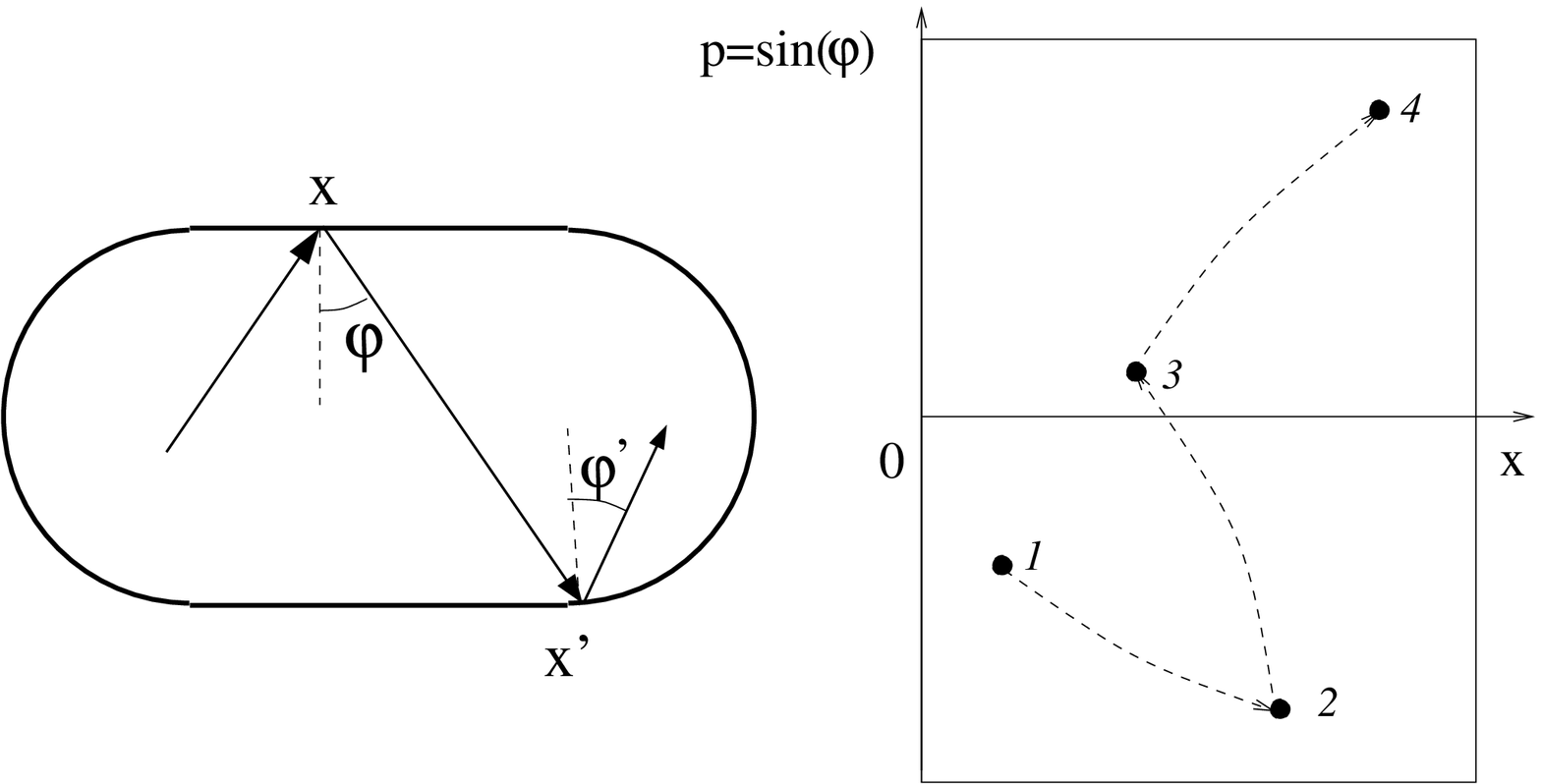}
\end{center}
\caption{Top: Poincar\'e section and the associated return map constructed
  from a Hamiltonian flow on $\cE_E$. Bottom: boundary map associated
  with the stadium billiard.}\label{f:Poincare}
\end{figure}
The advantage of studying maps instead of flows is multifold. 
Firstly, a map can be easily
constructed from a flow by considering a {\em Poincar\'e section}
$\Sigma$ transversal to the flow; the induced
return map $\kappa_\Sigma:\Sigma\to\Sigma$, together with the
return time, contain all the dynamical information on the
flow. Ergodic properties of chaotic maps are usually easier to
study than their flow counterpart. For billiards, the natural
Poincar\'e map to consider is the boundary map $\kappa_\Sigma$ defined on the phase space
associated with the boundary, $T^*\partial\Omega$. The ergodic
properties of this boundary map
were understood, and used to address the case of the billiard flow
itself \cite{BalMel08}.

\begin{figure}[htbp]
\begin{center}
\includegraphics[width=.8\textwidth]{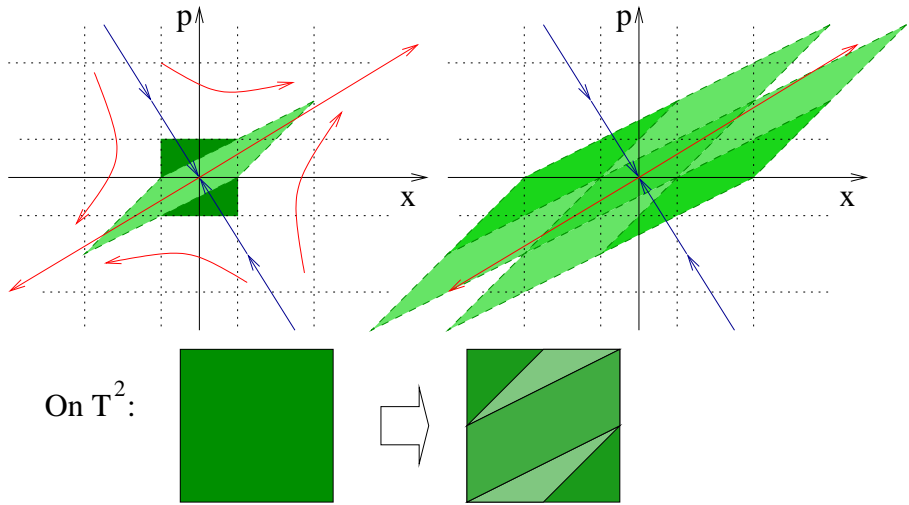}
\end{center}
\caption{Construction of Arnold's cat map $\kappa_{S_{cat}}$ on the
  2-torus, obtained by periodizing the linear transformation on
  $\IR^2$. The stable/unstable directions are shown (kindly
  provided by F.~Faure).}\label{f:catmap}
\end{figure}

Secondly, simple chaotic maps can be defined on low-dimensional phase
spaces, the most famous ones being the hyperbolic symplectomorphisms on
the 2-dimensional torus. These are defined by the action of a matrix
$S=\begin{pmatrix}a&b\\c&d\end{pmatrix}$ with integer entries,
determinant unity and trace $a+d>2$ (equivalently, $S$ is unimodular and
hyperbolic). Such a matrix obviously acts on
$\vx=(x,p)\in\t2$ linearly, through 
\be\label{e:catmap}
\kappa_S(\vx)=(ax+bp,cx+dp)\bmod 1\,.
\ee
A schematic view of $\kappa_S$ for the famous {\em Arnold's cat map}
$S_{cat}=\begin{pmatrix}1&1\\1&2\end{pmatrix}$ is displayed in Fig.~\ref{f:catmap}.
The hyperbolicity condition implies that
the eigenvalues of $S$ are of the form $\{e^{\pm\lambda}\}$ for some
$\lambda>0$. As a result, $\kappa_S$ has the Anosov property: at each point
$\vx$, the tangent space
$T_\vx\t2$ splits into stable and unstable subspaces, identified with the
eigenspaces of $S$, and $\pm\lambda$ are the Lyapunov exponents. 
Many dynamical properties of $\kappa_S$ can be explicitly computed. 
For instance, every rational point $\vx\in\t2$ is periodic, and the
number of periodic orbits of period $\leq n$ grows like $e^{\lambda
  n}$ (thus $\lambda$ also measures the complexity of the map).
This linearity also results in the fact that the decay of correlations (for smooth observables)
is superexponential.
To obtain a {\em generic} Anosov diffeomorphism of the 2-torus, one can
smoothly perturb the linear map $\kappa_S$: given a
Hamiltonian $H\in C^\infty(\t2)$, the composed map
$\Phi^\eps_H\circ\kappa_S$ remains Anosov
if $\eps$ is small enough, due to the structural stability of
Anosov diffeomorphisms. 
For any such Anosov diffeomorphism, the decay
of correlations is exponential.
\begin{figure}[htbp]
\begin{center}
\includegraphics[width=.7\textwidth]{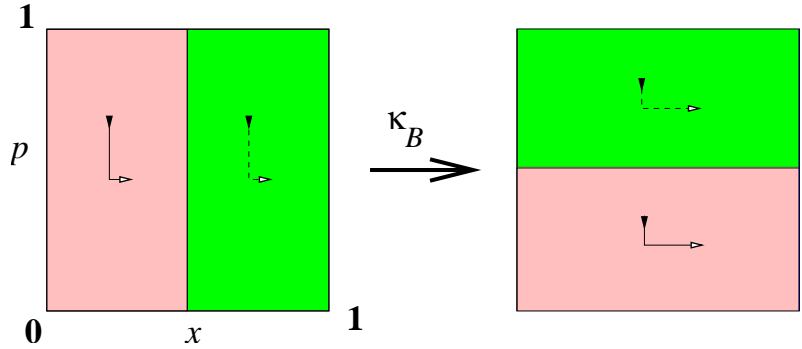}
\end{center}
\caption{Schematic view of the baker's map \eqref{e:baker}. The arrows
show the contraction/expansion directions.}\label{f:baker}
\end{figure}
Another family of canonical maps on the torus was also much investigated, namely the
so-called {\em baker's maps}, which are piecewise linear. The simplest
(symmetric) baker's map is defined by
\be\label{e:baker}
\kappa_B(x,p)=\begin{cases}(2x\bmod 1,\frac{p}{2}),&0\leq x<1/2,\\(2x\bmod 1,\frac{p+1}{2}),& 1/2\leq x<1.\end{cases}
\ee
This map is conjugate to a very simple
 {\em symbolic dynamics}, namely the shift on two
symbols. Consider the binary expansions of the
coordinates $x=0.\alpha_1\alpha_2\cdots$,
$p=0.\beta_1\beta_2\cdots$. Then the map $(x,p)\mapsto \kappa_B(x,p)$ is
equivalent with the shift to the left on the bi-infinite sequence
$\cdots\beta_2\beta_1\cdot\alpha_1\alpha_2\cdots$. This conjugacy allows to
easily identify all periodic
orbits, construct a large set of nontrivial invariant probability
measures, and prove that the map is ergodic and mixing w.r.t. these measures.
All trajectories not meeting the discontinuity lines are uniformly
hyperbolic. 

Simple canonical maps have also be defined on the 2-sphere phase space
(like the kicked top), but their chaotic properties have, to my knowledge,
not been rigorously proven. Their quantization has been
intensively investigated, mostly numerically \cite{Haake01}.

\subsubsection{Quantum maps on the 2-dimensional torus\label{s:quantum-maps}}

As opposed to the case of Hamiltonian flows, there is no natural
rule to quantize a canonical map on a compact phase space
$\cP$. Already, associating a quantum Hilbert space to this phase
space is not obvious. 
Therefore, from the very beginning, quantum maps have been defined through
somewhat arbitrary (or rather, {\em ad hoc}) procedures, often specific
to the considered map $\kappa:\cP\to\cP$. Still these
recipes are always required to satisfy a certain number of properties:
\begin{itemize}
\item one needs a {\em sequence} of Hilbert spaces
  $(\hn)_{N\in\IN}$ of dimensions $N$. Here $N$ is
  interpreted as the inverse of Planck's constant, in agreement
  with the heuristics that {\em each quantum state occupies a
  volume $\hbar^d$ in phase space}. We also want to quantize
observables $f\in C(\cP)$ into hermitian operators $\hat f_N$ on $\hn$.
\item  For each $N\geq 1$, the
quantization of $\kappa$ is given by a unitary propagator
$U_N(\kappa)$ acting on $\hn$. The whole family $(U_N(\kappa))_{N\geq
  1}$ is called the quantum map associated with $\kappa$. 
\item in the semiclassical limit
$N\sim\hbar^{-1}\to\infty$, this propagator satisfies some
form of quantum-classical correspondence. Namely, for some (large enough) family
of observables $f$ on $\cP$, we should have
\be\label{e:Egorov-map}
\forall n\in\IZ,\qquad U_N^{-n}\,\hat f_N\,U_N^n = \widehat{(f\circ\kappa^n)}_N
+\cO_n(N^{-1})\qquad\text{as $N\to\infty$.}
\ee
\end{itemize}
The condition \eqref{e:Egorov-map} is the analogue of the Egorov
property \eqref{e:Egorov} satisfied by the
propagator $U_\hbar$ associated with a quantum Hamiltonian, which
quantizes the stroboscopic map $\vx\mapsto \Phi^1_H(\vx)$.

Let us briefly summarize the explicit construction of the quantizations
$U_N(\kappa)$, for the maps $\kappa:\t2\to\t2$ presented in the
previous section. Let us start by constructing the quantum Hilbert
space.
One can see $\t2$ as the quotient of the phase space
$T^*\IR=\IR^2$ by the discrete tranlations 
$\vx\mapsto \vx+\vn,\quad \vn\in\IZ^2$.
Hence, it is natural to construct quantum states on $\t2$ by starting
from states $\psi\in \cS'(\IR)$, and requiring the following
periodicity properties
$$
\psi(x+n_1)=\psi(x),\quad (\cF_\hbar\psi)(p+n_2)=(\cF_\hbar\psi)(p),\quad n_1,n_2\in\IZ.
$$
It turns out that these two conditions can be satisfied only if
$\hbar=(2\pi N)^{-1}$, $N\in\IN$, and the corresponding distributions then
form a vector space $\hn$ of dimension $N$. A basis of this space is given by the Dirac combs
\be\label{e:e_j}
e_\ell(x)=\frac{1}{\sqrt{N}}\sum_{\nu\in\IZ}\delta(x-\frac{\ell}{N}-\nu),\qquad \ell=0,\ldots,N-1\,.
\ee
It is natural to equip $\hn$ with the hermitian structure for
which the basis $\{e_\ell,\,\ell=0,\ldots,N-1\}$ is orthonormal. The
components of a vector $\psi\in\cH_N$ in this ``position basis'' will be denoted by
$\psi(\ell/N)=\la e_\ell,\psi\ra$.

Let us now explain how the symplectomorpisms $\kappa_S$ are
quantized \cite{HanBer80}. Given a unimodular matrix $S$,
its action on $\IR^2$ can be generated by
a quadratic polynomial $H_S(x,p)$; this action can thus be quantized into the
unitary operator $U_\hbar(S)=\exp(-i\hat
H_{S,\hbar}/\hbar)$ on $L^2(\IR)$. This operator also
acts on distributions $\cS'(\IR)$, and in particular on the
finite subspace $\hn$. Provided the
matrix $S$ satisfies some ``checkerboard
condition'', one can show (using group theory) that the action of
$U_\hbar(S)$ on $\hn$ {\em preserves that space}, and acts on it
through a unitary matrix $U_N(S)$. The family
of matrices $(U_N(S))_{N\geq 1}$ defines the quantization
of the map $\kappa_S$ on $\t2$.
Group theory also implies that an {\em exact} quantum-classical correspondence
holds (that is, the remainder term in
\eqref{e:Egorov-map} vanishes), and has other important consequences regarding
the operators $U_N=U_N(\kappa_S)$ (for each $N$ the
matrix $U_N$ is periodic, of period $T_N\leq 2N$).
Explicit expressions for the matrices $U_N(S)$ can be
worked out, they depends sensitively on the factorization into primes of the integer $N$. 

\medskip

The construction of the quantized baker's map \eqref{e:baker}
proceeds differently. An Ansatz was proposed by
Balasz-Voros \cite{BalVor89}, which takes the following form in the
basis \eqref{e:e_j}\footnote{We assume that $N$ is an even integer}:
\be
U_N(\kappa_B)=F_N^*\begin{pmatrix}F_{N/2}&\\ &F_{N/2}\end{pmatrix}\,,
\ee
where $F_N$ is the $N$-dimensional discrete Fourier transform. This
Ansatz is obviously unitary. It was obtained by discretizing the phases $e^{i\varphi(p',x)/\hbar}$,
where the function $\varphi(p',x)$ locally generates\footnote{This means
  that, in some region, $(x',p')=\kappa_B(x,p)$ is defined by solving
  $x'=\partial_{p'}\varphi(p',x)$, $p=\partial_x\varphi(p',x)$.} the map \eqref{e:baker}.
A proof that the matrices $U_N(\kappa_B)$
satisfy the Egorov property \eqref{e:Egorov-map} was given in \cite{DENW06}. 

Once we have constructed the matrices $U_N(\kappa)$
associated with a chaotic map $\kappa$, their eigenstates
$\{\psi_{N,j},\,j=1,\ldots,N\}$ enjoy the r\^ole of quantum chaotic
eigenstates. They are of quite different nature from the eigenstates of the
Laplacian on a manifold or a billiard: while the latter belong
to $L^2(M)$ or $L^2(\Omega)$ (and are actually smooth functions), the
eigenstates $\psi_{N,j}$ are $N$-dimensional vectors. 
Still, part of ``quantum chaos'' has consisted in
developing common tools to analyze these eigenstates, in spite of the different functional settings.

\section{Macroscopic description of the eigenstates\label{s:macro}}

In this section we study the {\em macroscopic} localization properties of chaotic
eigenstates. Most of the results are mathematically rigorous. 
In the case of the semiclassical Laplacian \eqref{e:Laplacian} on a billiard $\Omega$ we ask
the following question:
\begin{quote}
Consider a sequence
$(\psi_\hbar)_{\hbar\to 0}$ of normalized eigenstates of $\hat H_\hbar$, with
energies $E_\hbar\approx 1/2$. For $A\subset\Omega$ a fixed subdomain, 
what is the probability that the particle described by the stationary state
$\psi_\hbar$ lies inside $A$? How do the probability weights
$\int_{A}|\psi_\hbar(x)|^2\,dx$ behave when $\hbar\to 0$?
\end{quote}
This question is quite natural, when contemplating eigenstate
plots like in Fig.~\ref{f:stadium-eigen}. Here, by {\em macroscopic} we mean
that the domain $A$ is kept fixed while $\hbar\to 0$.

One can obviously generalize the question to integrals of the type
$\int_{\Omega}f(x)\,|\psi_\hbar(x)|^2\,dx$, with $f(x)$ a continuous
test function on $\Omega$. Assuming the volume of $\Omega$ is
normalized, this integral can be interpreted as the
{\em quantum average}\footnote{These are also referred to as the
  (diagonal) {\em matrix elements} for the ``matrix'' $\hat f_\hbar$.} 
$\la\psi_\hbar,\hat f_\hbar\,\psi_\hbar\ra$, where the
quantum observable $\hat f_\hbar$ is just the multiplication operator
by $f(x)$. It proves useful to extend the question to
phase space observables $f(x,p)$\footnote{If $\Omega\subset\IR^d$ is replaced by a compact riemannian
manifold $M$, we use some quantization scheme $f\mapsto \hat f_\hbar$.}: 
what is the behaviour of the quantum averages
\be
\mu^W_{\psi_\hbar}(f)\defeq\la \psi_\hbar,\hat f_\hbar\psi_\hbar\ra,\qquad f\in
C^\infty(T^*\Omega)\,,\quad\text{in the limit $\hbar\to 0$?}
\ee
Since the quantization procedure $f\mapsto \hat f_\hbar$ is linear, 
these averages define a distribution $\mu^W_{\psi_\hbar}=W_{\psi_\hbar}(\vx)\,d\vx$ on
$T^*\Omega$, called the Wigner distribution
of the state $\psi_\hbar$  (the density $W_{\psi_\hbar}(\vx)$ is called the
Wigner function). Although this distribution is generally not
positive, it is interpreted as a quasi-probability density
describing the state $\psi_\hbar$ in {\em phase space}.

On the euclidean space, one can define a {\em nonnegative} phase space
density associated to the state
$\psi_\hbar$: the {\em Husimi} measure (and function):
\be\label{e:Husimi}
\mu^H_{\psi_\hbar}=\cH_{\psi_\hbar}(\vx)\,d\vx,\qquad 
\cH_{\psi_\hbar}(\vx) = (2\pi\hbar)^{-d/2}\,|\la\varphi_{\bx},\psi_\hbar\ra|^2\,,
\ee
where the Gaussian
wavepackets $\varphi_{\vx}\in L^2(\IR^d)$ are defined in \eqref{e:wavepacket}. This
measure can also be obtained by convolution of $\mu^W_{\psi_\hbar}$
with the Gaussian kernel $e^{-|\vx-\vy|^2/\hbar}$. 
In these notes our phase space plots show Husimi functions.

These questions lead us to the notion of phase space localization, or {\em
  microlocalization}\footnote{The prefix
  {\em micro} mustn't mislead us: we are still dealing with {\em
    macroscopic} localization properties of $\psi_\hbar$!}. We will
say that the family of states
$(\psi_\hbar)$ is microlocalized inside a set $B\subset
T^*\Omega$ if, for any smooth observable $f(x,p)$ vanishing near $B$, the
quantum averages $\la \psi_\hbar,\hat f_\hbar\psi_\hbar\ra$
decrease faster than any power of $\hbar$ when $\hbar\to 0$. 

Microlocal properties are not easy to guess from plots of the spatial
density $|\psi_j(x)|^2$ like in
Fig.~\ref{f:stadium-eigen}, but they are natural to study if we want
to connect quantum to classical mechanics, since the latter takes
place in phase space rather than in position space. Indeed, the major
tool we will use is the
quantum--classical correspondence \eqref{e:Egorov}; for all the flows we consider,
a generic {\em spatial} test function $f(x)$ evolves through $\Phi^t$
into a genuine {\it phase
space} function $f_t(x,p)$. Microlocal properties are easier to visualize on
2-dimensional phase spaces (see below the figures on the 2-torus).

\subsection{The case of completely integrable systems}
In order to motivate our further discussion of chaotic eigenstates, let us
first recall a few facts about the antipodal systems, namely
completely integrable Hamiltonian flows. For such systems, the energy shell $\cE_E$ is foliated
by $d$-dimensional invariant Lagrangian tori. Each such torus is characterized
by the values of $d$ independent {\em invariant actions}
$I_1,\ldots,I_d$, so let us call such a torus $T_{\vec I}$. 
The WKB theory allows one to explicitly construct, in the
semiclassical limit, precise
{\em quasimodes} of $\hat H_\hbar$ associated with some of these tori\footnote{The
  ``quantizable'' tori $T_{\vec I}$ satisfy
  Bohr-Sommerfeld conditions $I_i=2\pi\hbar (n_i + \alpha_i)$, with
  $n_i\in\IZ$ arbitrary and $\alpha_i\in [0,1]$ fixed Maslov indices.},
that is normalized states $\psi_{\vec I}=\psi_{\hbar,\vec I}$ satisfying
\be
\hat H_\hbar \,\psi_{\vec I} = E_{\vec I}\, \psi_{\vec I} + \cO(\hbar^\infty),
\ee
with energies $E_{\vec I}\approx E$. Such a Lagrangian (or WKB)
state $\psi_{\vec I}$ takes the following form (away from caustics):
\be\label{e:Lagrangian}
\psi_{\vec I}(x)=\sum_{\ell=1}^L A_\ell(x;\hbar)\,\exp(i S_\ell(x)/\hbar)\,.
\ee
Here the functions $S_\ell(x)$ are (local) generating
functions\footnote{Above some neighbourhood $U\in\IR^d$, the
torus $T_{\vec I}$ is the union of $L$ lagrangian leaves
$\{(x,\nabla S_\ell(x)),\,x\in U\}$, $\ell=1,\ldots,L$} for $T_{\vec I}$, and each
$A_\ell(x;\hbar)=A_\ell^0(x)+\hbar A_\ell^1(x)+\cdots$ is a
smooth amplitude.

From this explicit expression one can check that the state $\psi_{\vec I}$ is
microlocalized on $T_{\vec I}$. On the other hand, our knowledge
of $\psi_\hbar$ is much more precise than the latter fact; indeed, one can
easily construct states microlocalized on $T_{\vec I}$, which are very
different from the Lagrangian states $\psi_{\vec I}$. For instance, for
any point $\vx_0=(x_0,p_0)\in T_{\vec I}$ the Gaussian wavepacket (or
coherent state)
\be\label{e:wavepacket}
\varphi_{\vx_0}(x)=(\pi\hbar)^{-1/4}\,e^{-|x-x_0|^2/2\hbar}\,e^{ip_0\cdot x/\hbar}
\ee
is microlocalized on the single point $\vx_0$, and therefore
on $T_{\vec I}$. This example just reflects the fact that a statement
about microlocalization of a sequence of states provides much less information
than a formula like \eqref{e:Lagrangian}. 
\begin{figure}[htbp]
\begin{center}
\includegraphics[angle=-90,width=0.35\textwidth]{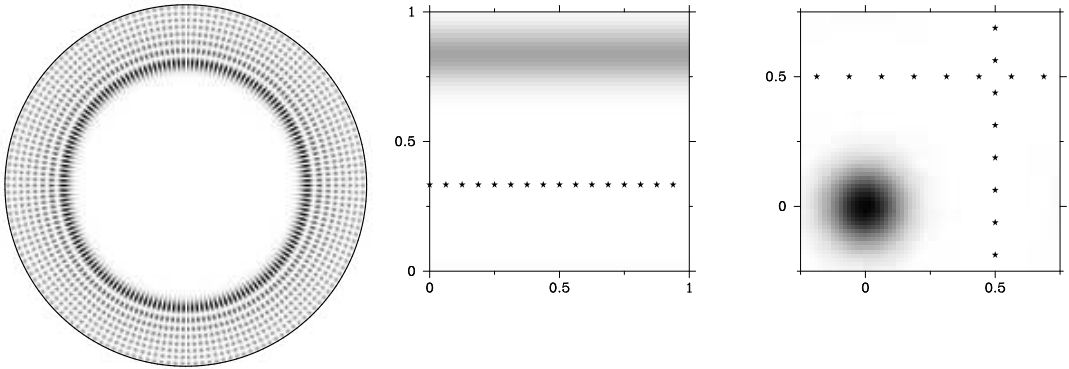}
\includegraphics[angle=-90,width=0.6\textwidth]{lagrang+coherent.ps}
\caption{Left: one eigenmode of the circle billiard, microlocalized on
a torus $T_{\vec I}$ (the invariants $I_1,I_2$ can be taken to be the energy and the
angular momentum). Center: Husimi density of a simple Lagrangian state
on $\t2$, a momentum eigenstate microlocalized on the
lagrangian $\{\xi=\xi_0\}$. 
Right: Husimi density of the Gaussian wavepacket \eqref{e:wavepacket}. Stars
denote the zeros of the Husimi density. (reprinted from \cite{NonVor98})}\label{f:integ}
\end{center}
\end{figure}

\subsection{Quantum ergodicity}

In the case of a fully chaotic system, there is no explicit formula
describing the eigenstates, or even the quasimodes, of $\hat H_\hbar$. However,
macroscopic informations on the eigenstates can be obtained
indirectly, using the quantum-classical correspondence. The main result on this
question is a quantum analogue of the ergodicity property \eqref{e:ergodicity} of the
classical flow, and it is a consequence of this property. For this
reason, it has been named {\em Quantum Ergodicity} by
Zelditch. Loosely speaking, this property states that {\em almost} all
eigenstates $\psi_\hbar$ with $E_\hbar\approx E$ will become {\em
  equidistributed} on the energy
shell $\cE_E$, in the semiclassical limit, provided the classical flow on
$\cE_E$ is ergodic. We give below the version of
the theorem in the case of the Laplacian on a compact
riemannian manifold, using the notations of \eqref{e:Helmholtz}.
\begin{Theorem}\textbf{[Quantum ergodicity]}\label{thm:QE}
Assume the geodesic flow on $(M,g)$ is ergodic w.r.to the
Liouville measure. Then, for any orthonormal eigenbasis
$(\psi_j)_{j\geq 0}$ of the Laplacian, there exists a subsequence
$\cS\subset\IN$ of density $1$ (that is,
$\lim_{J\to\infty}\frac{\#(\cS\cap [1,J])}{J}=1$), such that 
$$
\forall f\in C^\infty_c(M),\qquad \lim_{j\in \cS,j\to\infty} 
\la \psi_j,\hat f_{\hbar_j}\psi_j\ra = \int_{\cE}f(x,p)\,d\mu_L(x,p)\,,
$$
where $\mu_L$ is the normalized Liouville measure on $\cE=\cE_{1/2}$,
and $\hbar_j=k_j^{-1}$ the inverse wavevector of $\psi_j$.
\end{Theorem}
The statement of this theorem was first given by Schnirelman (using
test functions $f(x)$) \cite{Schnirel74}, the complete proof was
obtained by Zelditch in the case of manifolds of constant negative curvature
$\Gamma\backslash\IH$ \cite{Zel87}, and the general case was then
proved by Colin de Verdi\`ere \cite{CdV85}.
This theorem is ``robust'': it has been extended
to
\begin{itemize}
\item quantum ergodic billiards \cite{GL93,ZZ96}
\item quantum Hamiltonians $\hat H_\hbar$, such that the flow $\Phi^t_H$ is
ergodic on $\cE_E$ in some energy interval \cite{HMR87}
\item quantized ergodic diffeomorphisms on the torus \cite{BDB96} or
  on more general compact phase spaces \cite{Zel97}
\item a general framework of $C^*$ dynamical systems \cite{Zel-C*96}
\item a family of quantized ergodic maps with discontinuities
  \cite{MOK05}, in particular the
  baker's map \cite{DENW06}
\item certain quantum graphs \cite{BerKeatSmil07}
\end{itemize}
Let us sketch the proof of the above theorem.
We want to study the statistical distribution of the matrix elements
$\mu^W_{j}(f)=\la \psi_j,\hat f_{\hbar_j}\psi_j\ra$ in the range
$\{k_j\leq K\}$, with $K\gg 1$. The first step is to estimate the
average of this distribution. It is estimated by the generalized
Weyl law:
\be\label{e:general-Weyl}
\sum_{k_j\leq K}\mu_{j}(f)\sim \frac{\Vol(M)\sigma_d}{(2\pi)^d}\,K^{d}\,\int_\cE f\,d\mu_L,\quad K\to\infty\,,
\ee
where $\sigma_d$ is the volume of the unit ball in $\IR^d$ \cite{Horm68}. 
In particular, this asymptotics allows to count the number of
eigenvalues $k_j\leq K$:
\be\label{e:Weyl}
\#\{k_j\leq K\}\sim \frac{\Vol(M)\sigma_d}{(2\pi)^d}\,K^d,\quad K\to\infty,
\ee
and shows that the average of the distribution $\{\mu^W_j(f),\,k_j\leq K\}$ converges
to the phase space average $\mu_L(f)$ when $K\to\infty$.

Now, we want to show that the distribution is
concentrated around its average. This can be done by estimating its variance
$$
\Var_K(f)\defeq
\frac{1}{\#\{k_j\leq K\}}\sum_{k_j\leq K}|\la \psi_j,(\hat f_{\hbar_j}-\mu_L(f))\psi_j\ra|^2\,,
$$
which is often called the {\em quantum variance}.
Because the $\psi_j$ are eigenstates of the propagator $U^t_\hbar$, we may replace $\hat
f_{\hbar_j}$ by its time average up to some large time $T$,
$$
\hat f_{T,\hbar_j}=
\frac{1}{2T}\int_{-T}^T U_{\hbar_j}^{-t}\hat f_{\hbar_j}U_{\hbar_j}^{t}\,dt\,,
$$
without modifying the matrix elements. Then we apply the
Cauchy-Schwarz inequality:
$$
|\la\psi_j,A\psi_j|^2\leq \la\psi_j,A^*\,A\psi\ra,\quad \text{for any
  bounded operator $A$},
$$
to get the following upper bound for the variance:
$$
\Var_K(f)\leq \frac{1}{\#\{k_j\leq K\}}\sum_{k_j\leq K}\la \psi_j,(\hat f_{T,\hbar_j}-\mu_L(f))^*(\hat f_{T,\hbar_j}-\mu_L(f))\psi_j\ra\,.
$$
The Egorov theorem \eqref{e:Egorov} shows that the
product operator on the right
hand side is approximatly equal to the quantization of the function
$|f_{T}-\mu_L(f)|^2$, where $f_T$ is the time average of the classical observable
$f$. 
Applying the generalized Weyl law \eqref{e:general-Weyl} to this function, we get
$$
\Var_K(f)\leq \mu_L(|f_{T}-\mu_L(f)|^2)+\cO_T(K^{-1})\,,\quad K\to\infty.
$$
Finally, the ergodicity of the classical flow implies that
$\mu_L(|f_{T}-\mu_L(f)|^2)$ converges to zero when $T\to\infty$. By
taking $T$, and then $K$, large enough, the above right hand
side can be made arbitrary small. This proves that the quantum variance $\Var_K(f)$
converges to zero when $K\to\infty$. A standard Chebychev argument is then used to extract a
dense subsequence such that $\la\psi_\hbar,\hat f_\hbar \psi_\hbar\ra
\to \mu_L(f)$, and a final diagonal extraction procedure yields a
dense subsequence enjoying the same property for all smooth observables. $\hfill\square$

A more detailed discussion on the quantum variance and the distribution of matrix elements
$\{\mu_j(f),\,k_j\leq K\}$ will be given in \S\ref{s:distrib}.

\subsection{Beyond QE: Quantum Unique Ergodicity vs. strong  scarring}

\subsubsection{Semiclassical measures\label{s:mu_sc}}
Quantum ergodicity can be conveniently expressed by using the
concept of {\em semiclassical measure}. Remember that, using a duality
argument, we associate to each eigenstates $\psi_j$ a Wigner
distribution $\mu^W_j$ on phase space. The quantum ergodicity theorem~\ref{thm:QE} can be rephrased as follows:
\begin{quote}
There exists a density-$1$ subsequence $\cS\subset\IN$, such that the
sequences of Wigner distributions $(\mu^W_{j})_{j\in \cS}$ weak-$*$ (or vaguely) converges to
the Liouville measure on $\cE$.
\end{quote}
For any compact riemannian manifold $(M,g)$, the
sequence of Wigner distributions $(\mu^W_j)_{j\in\IN}$ remains in a
compact set in the weak-$*$ topology, so it is always possible to
extract an infinite subsequence $(\mu_j)_{j\in \cS}$ vaguely converging to a limit
distribution $\mu_{sc}$, that is 
$$
\forall f\in C^\infty_c(T^*M),\quad \lim_{j\in \cS,j\to\infty}\int
f\,d\mu^W_j = \int f\,d\mu_{sc}\,.
$$ 
Such a limit distribution is necessarily a probability measure
on $\cE$, and is called a {\em semiclassical measure} of the manifold
$M$. The quantum-classical correspondence implies that
$\mu_{sc}$ is {\em invariant} through the geodesic flow: 
$(\Phi^t)^*\mu_{sc}=\mu_{sc}$.
The semiclassical measure $\mu_{sc}$ represents the asymptotic
(macroscopic) phase space distribution of the eigenstates
$(\psi_j)_{j\in \cS}$. It is the major tool used in the mathematical
literature on chaotic eigenstates (see below). The definition can be
obviously generalized to any quantized Hamiltonian flow or canonical
map.

\subsubsection{Quantum Unique Ergodicity conjecture}
The quantum ergodicity theorem provides an incomplete information,
which leads to the following question:
\begin{Question}
Do {\em all} eigenstates become equidistributed in the semiclassical
limit? 
Equivalently, is the Liouville measure the {\em unique} semiclassical
measure for the manifold $M$? On the opposite, are there {\em exceptional
subsequences} converging to a semiclassical measure $\mu_{sc}\neq \mu_L$?
\end{Question}
This question makes sense if the geodesic flow admits
invariant measures different from $\mu_L$ (that is, the flow is not
uniquely ergodic). Our central example, manifolds of negative
curvature, admit many different invariant measures, e.g. the singular
measures $\mu_\gamma$ supported on each of the (countably many)
periodic geodesics $\gamma$.

This question was already raised in \cite{CdV85}, where the author
conjectured that no subsequence of eigenstates can concentrate along
a single periodic geodesic, in other words  $\mu_{\gamma}$ cannot be a semiclassical measure.
Such an unlikely subsequence was later called a {\em strong scar} by Rudnick and
Sarnak \cite{RudSar94}, in reference to the {\em scars} discovered by
Heller on the stadium billiard (see \S\ref{s:scars}).
In the same paper the authors formulated a stronger conjecture:
\begin{Conj}\label{c:QUE}\textbf{[Quantum unique ergodicity]}
Let $(M,g)$ be a compact riemannian manifold with negative sectional
curvature. Then all high-frequency eigenstates of the Laplacian become equidistributed
with respect to the Liouville measure. Equivalently, the latter is the
unique semiclassical measure of $M$.
\end{Conj}
The term {\em quantum unique ergodicity} refers to the notion of
{\em unique ergodicity} in ergodic theory: a system is uniquely
ergodic system if it admits unique invariant probability measure. The
geodesic flow we are considering admits many invariant
measures, but the conjecture states that the corresponding quantum
system selects only one of them.

\subsubsection{Arithmetic quantum unique ergodicity\label{s:arithmetic}}
The QUE conjecture was motivated by the following
result obtained in the cited article. Let us recall that any invariant measure can be decomposed
into a convex combination of {\em ergodic} invariant measures:
\be\label{e:ergod-decompo}
\mu = \int_{\mathrm{Erg}}\mu_e\,d\nu(e)\,,
\ee
where $\{\mu_e,\ e\in \mathrm{Erg}\} $ span the set of ergodic probability measure, and $\nu$ is a
probability measure over this set.
The authors specifically considered
arithmetic surfaces, obtained by
quotienting the Poincar\'e disk $\IH$ by certain congruent co-compact
lattices $\Gamma$. As explained in \S\ref{s:geodesic}, on such a surface
it is natural to consider Hecke eigenstates, which are joint eigenstates of the Laplacian and the
(countably many) Hecke operators\footnote{The spectrum of the
  Laplacian on such a surface is believed to be simple; if this is the case, then an
  eigenstate of $\Delta$ is automatically a Hecke eigenstate.}.
It was shown in \cite{RudSar94} that for any semiclassical measure
$\mu_{sc}$ associated with a sequence of Hecke eigenstates, the
ergodic decomposition \eqref{e:ergod-decompo} of $\mu_{sc}$ does not charge any
ergodic component $\mu_{\gamma}$ associated with a periodic geodesic.
The methods of \cite{RudSar94} were refined by Bourgain and
Lindenstrauss \cite{BLi03}, who showed that 
the measure $\mu_{sc}$ of an $\eps$-thin tube along any geodesic segment is bounded
from above by $C\,\eps^{2/9}$; this bound implies
that the {\em Kolmogorov-Sinai entropy}\footnote{The notion of KS
  entropy will be explained in \S\ref{s:entropy}.} of almost every ergodic
component of $\mu_{sc}$ is bounded from below by $2/9$. 
Finally, using advanced ergodic theory methods, Lindenstrauss
completed the proof of the QUE conjecture in the arithmetic
context. 
\begin{Theorem}{\bf [Arithmetic QUE]}\cite{Linden06}
Let $M=\Gamma\backslash\IH$ be an arithmetic\footnote{For the precise
  definition of these surfaces, see \cite{Linden06}.} surface of constant
negative curvature. Consider an eigenbasis $(\psi_j)_{j\in\IN}$ of
Hecke eigenstates of $\Delta_M$. Then, the only semiclassical measure
associated with this sequence is the Liouville measure.
\end{Theorem}
Lindenstrauss and Brooks recently improved this
theorem \cite{BroLin10b,BroLin11}: the QUE result holds true, assuming the
$(\psi_j)$ are joint eigenstates of the Laplacian and of a (nonelementary) single Hecke
operator. Their proof uses a new delocalization estimate for
regular graphs \cite{BroLin10}, which yields the same positive entropy results
as in \cite{BLi03}.

\medskip

\begin{figure}[htbp]
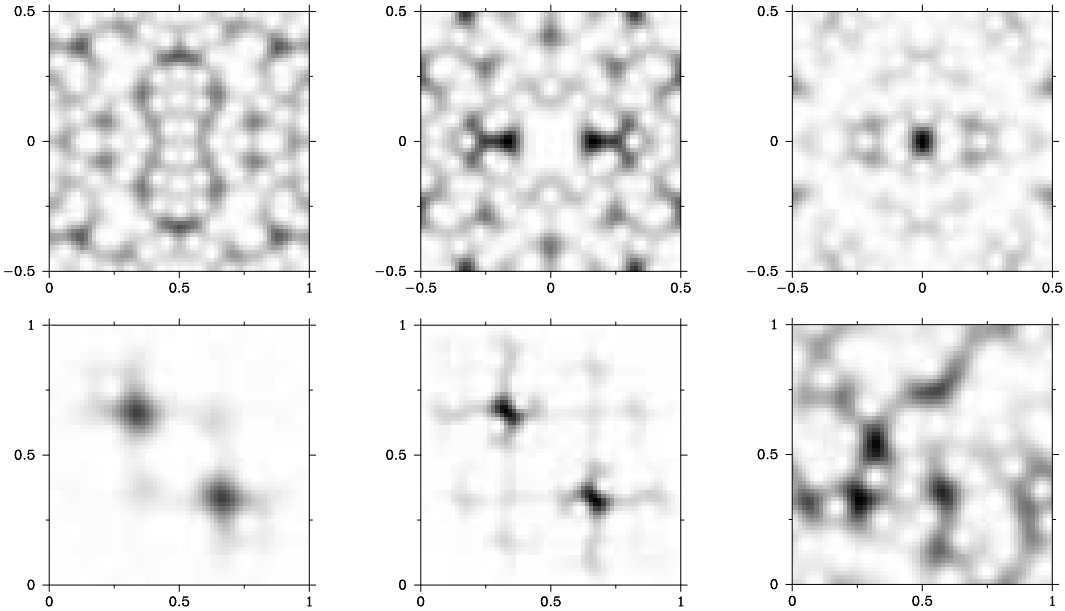

\begin{center}
\includegraphics[angle=-90,width=.95\textwidth]{cat-3states.ps}
\includegraphics[angle=-90,width=.95\textwidth]{baker-2states-random.ps}
\caption{Husimi functions of various states in $\hn$ (large values=dark regions). 
Top: 3 (Hecke) eigenstates of the quantum cat map $U_N(S_{DEGI})$, for
    $N=107$. Bottom left, center: 2 eigenstates of the quantum baker $U_N(\kappa_B)$ scarred on
    the period-2 orbit, for $N=48$ and $N=128$. Bottom right: random
    state \eqref{e:random-t2} for $N=56$. 
 (Reprinted from \cite{NonVor98})\label{f:cat-Hecke}}
\end{center}
\end{figure}

This arithmetic QUE result was preceded by a similar statement for the
hyperbolic symplectomorphisms on $\t2$ introduced in
\S\ref{s:quantum-maps}. These quantum maps $U_N(S_0)$ have the nongeneric
property to be periodic, so one has an explicit expression for
their eigenstates. It was shown in \cite{DEGI95} that for
a certain family hyperbolic matrices $S_0$ and certain sequences of
{\em prime} values of $N$, all the eigenstates of $U_N(S_0)$ become
equidistributed in the semiclassical limit, with an explicit bound on
the rate of equidistribution; they used explicit expressions for the
eigenstates, in terms of certain exponential sums. Some eigenstates, corresponding to the matrix
$S_{DEGI}=\begin{pmatrix}2&1\\3&2\end{pmatrix}$ 
are plotted in Fig.~\ref{f:cat-Hecke}, using the Husimi representation
on $\t2$.
A few years later, Kurlberg and Rudnick \cite{KurRud01}showed that this equidistribution
holds true provided the period of the propagator $U_N(S_0)$ is large
enough\footnote{basically, the period needs to be larger than
  $N^{1/2+\eps}$ for some $\eps>0$.}, which happens to be the case for almost all
  integers $N$ in the limit $N\to\infty$. 

A little earlier \cite{KurRud00}, they had constructed, attached to a given symplectomorphism
$S_0$ and any value of $N$, a
finite commutative family of operators $\{U_N(S') S'\in
\cC(S_0,N)\}$ including $U_N(S_0)$, which they called ``Hecke
operators'' by analogy with the case of arithmetic surfaces. They then
considered specifically the joint (``Hecke'') eigenbases of this family, and
proved QUE in this ``arithmetic setting'':
\begin{Theorem}\cite{KurRud00}
Let $S_0\in SL_2(\IZ)$ be a quantizable symplectic matrix. For each
$N>0$, consider a Hecke eigenbasis $(\psi_{N,j})_{j=1,\ldots,N}$ of
the quantum map $U_N(S_0)$. Then, for any observable $f\in
C^\infty(\t2)$ and any $\eps>0$, we have
$$
\la \psi_{N,j},\hat f_{N}\psi_{N,j}\ra = \int f\,d\mu_L + \cO_{f,\eps}(N^{-1/4+\eps}),
$$
where $\mu_L$ is the Liouville (or Lebesgue) measure on $\t2$.
\end{Theorem}
The eigenstates considered in \cite{DEGI95} were already particular instances of
these Hecke eigenstates.

In view of these positive results, it is tempting to generalize
the QUE conjecture to other chaotic systems.
\begin{Conj}{\bf [Generalized QUE]}\label{e:QUE-general}
Let $\Phi^t_H$ be an ergodic Hamiltonian flow on some energy shell
$\cE_E$. Then, all eigenstates $\psi_{\hbar,j}$ of $\hat H_\hbar$ of
energies $E_{\hbar,j}\approx E$ become equidistributed when $\hbar\to 0$.

Let $\kappa$ be a canonical ergodic map on $\t2$, and $(U_N(\kappa))_{N\in\IN}$ an
associated quantum map. Then, all the
eigenstates of $U_N(\kappa)$ become equidistributed on $\t2$
when $N\to\infty$.
\end{Conj}
An intensive numerical study for eigenstates of a Sinai-like billiard was carried on by
Barnett \cite{Bar06}. It seems to confirm QUE for this system.

In the next subsection we will exhibit particular systems for which this conjecture {\bf fails}.

\subsubsection{Counterexamples to QUE for quantum maps\label{s:counter}}

In this section we will exhibit sequences of eigenstates of certain
quantized chaotic maps, converging to
semiclassical measures different from $\mu_L$, thus disproving the
above conjecture.

Let us continue our discussion of symplectomorphisms on $\t2$.
We recall that for any $N\geq 1$, the quantized symplectomorphism
$U_N(S_0)$ is
periodic (up to a global phase) of period $T_N\leq 3N$, so that its eigenvalues are
essentially $T_N$-roots of unity. For values of $N$ such that $T_N\ll
N$, the spectrum of $U_N(S_0)$ is very degenerate, in which
case imposing the eigenstates to be of Hecke type becomes a strong
requirement. 
We already mentioned that, provided the
period is not too small (namely, $T_N\gg N^{1/2+\eps}$, which is
the case for almost all values of $N$), then QUE holds for any
eigenbasis \cite{KurRud01}.
On the opposite, there exist (sparse)
values of $N$, for which the period can be as small as
$T_N\sim C\,\log N$, so that the eigenspaces have huge
dimensions $\sim C^{-1}N/\log N$. This freedom allowed Faure,
De~Bi\`evre and the author to explicitly
construct eigenstates with different localization
properties \cite{FNdB03,FN04}.
\begin{Theorem}
Take $S_0\in SL_2(\IZ)$ a (quantizable) hyperbolic matrix. Then, there
exists an infinite (sparse) sequence $\cS\subset\IN$ such that, for
any periodic orbit $\gamma$ of $\kappa_{S_0}$, one can construct a
sequence of eigenstates $(\psi_N)_{N\in \cS}$ of $U_N(S_0)$ associated
with the semiclassical measure 
\be
\label{e:half-loc}\mu_{sc}=\frac12\mu_\gamma+\frac12\mu_L\,.
\ee
More generally, for any $\kappa_{S_0}$-invariant measure $\mu_{inv}$,
one can construct sequences of eigenstates associated
with the semiclassical measure 
$$
\mu_{sc}=\frac12\mu_{inv}+\frac12\mu_L.
$$
\end{Theorem}
This result provided the first counterexample to the generalized QUE
conjecture. The eigenstates ``converging to''
$\frac12\mu_\gamma+\frac12\mu_L$ can be called
{\em half-scarred}. The coefficient $1/2$ in front of the singular
component of $\mu_\gamma$ was shown to be optimal \cite{FN04}, a
phenomenon which was then better understood when considering the
entropy of $\mu_{sc}$ (see
\S\ref{s:entropy}). 

Let us briefly explain the construction of eigenstates half-scarred on
a fixed point $\gamma=\{\vx_0\}$ of $\kappa_{S_0}$. They are obtained by projecting on any
eigenspace the Gaussian wavepacket $\varphi_{\vx_0}$ (see
\eqref{e:wavepacket}). 
Each spectral projection can be expressed as a
linear combination of the evolved states $U_N(S_0)^n\varphi_{\vx_0}$, for
$n\in [-T_N/2,T_N/2-1]$. Now, we use the fact that, for $N$ in an
infinite subsequence $\cS\subset\IN$, the
period $T_N$ of the operator $U_N(S_0)$ is approximately the double of
the {\em Ehrenfest time}
\be\label{e:Ehrenfest}
T_E=\frac{\log\hbar^{-1}}{\lambda}\,,
\ee
(here $\lambda$ is the positive Lyapunov exponent).
The above linear combination can be split in two components:
during the time range
$n\in[-T_E/2,T_E/2]$ the states $U_N(S_0)^n\varphi_{\vx_0}$ remain
microlocalized at the origin; on the opposite, for times
$T_E/2<|n|\leq T_E$, these
states expand along long stretches of stable/unstable manifolds, and densely fill the
torus. As a result, the sum of these two components is half-localized, half-equidistributed.
$\hfill\square$

\begin{figure}
\begin{center}
\includegraphics[angle=-90,width=0.3\textwidth]{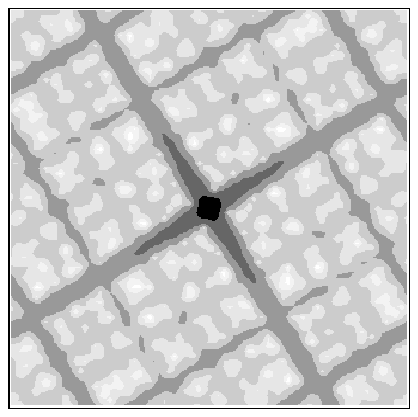}\hspace{.5cm}
\includegraphics[angle=-90,width=0.3\textwidth]{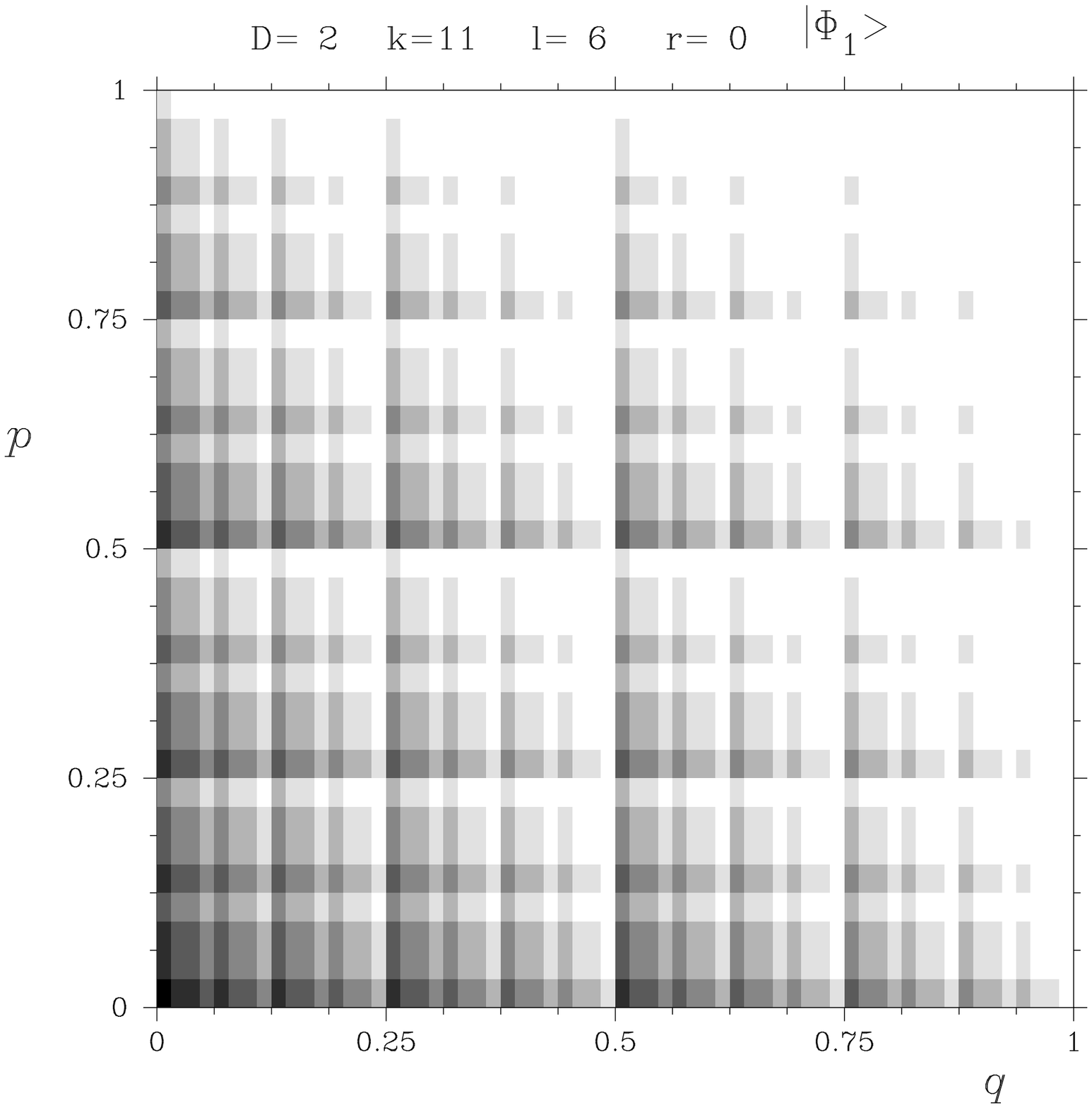}\hspace{.5cm}
\includegraphics[angle=-90,width=0.3\textwidth]{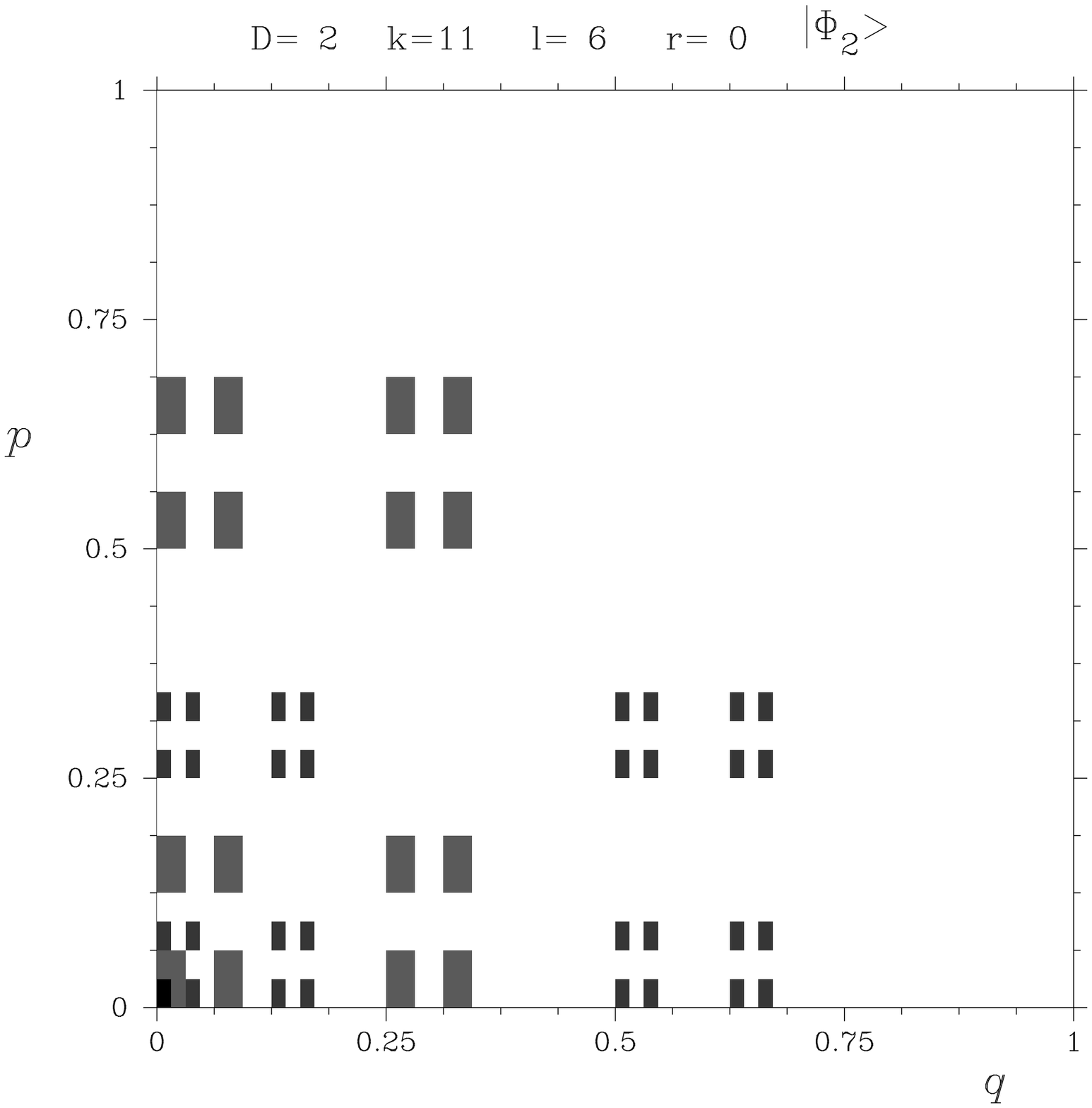}
\caption{Left: Husimi density of an eigenstate of
  $U_N(S_{cat})$, half-scarred on a fixed point. Notice the hyperbolic
structure around the fixed point  (Repr. from \cite{FNdB03}). Center, right: two eigenstates of the
  Walsh-quantized baker's map (we plot their ``Walsh-Husimi
  density''); the corresponding semiclassical
  measures are purely fractal (Repr. from \cite{AnaNo07-1}).\label{f:strong-scars}}
\end{center}
\end{figure}

In Fig.~\ref{f:strong-scars} (left) we plot the Husimi density
associated with one half-localized eigenstate of the quantum cat map
$U_N(S_{cat})$.

\medskip

A nonstandard (Walsh-) quantization of the 2-baker's map was
constructed in \cite{AnaNo07-1}, with properties similar to the above
quantum cat map. It allows to exhibit semiclassical measures $\frac12\mu_\gamma+\frac12\mu_L$ as in the above
case, but also purely fractal
semiclassical measures void of any Liouville component (see Fig.~\ref{f:strong-scars}).

Studying hyperbolic toral symplectomorphisms on $\IT^{2d}$ for $d\geq 2$,
Kelmer \cite{Kelmer06} identified eigenstates microlocalized on a proper
subtorus of dimension $\geq d$. He extended his analysis to certain nonlinear perturbations
of $\kappa_{S_0}$.
Other very explicit counterexamples to generalized QUE were constructed in \cite{CKST08}, based on
interval-exchange maps of the interval (such maps are ergodic, but have zero Lyapunov exponents).

\subsubsection{Counterexamples to QUE: bouncing-ball modes of the stadium billiard}

\begin{figure}[ht]
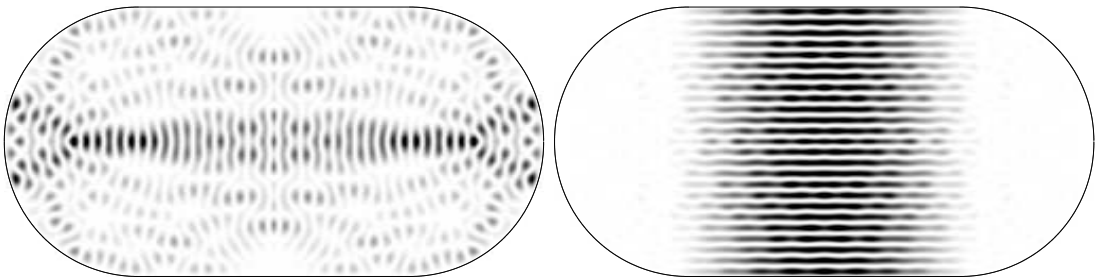

\begin{center}
\includegraphics[angle=-90,width=.48\textwidth]{stadium39ee-scar.ps}
\includegraphics[angle=-90,width=.48\textwidth]{stadium39ee-bb.ps}
\caption{Two eigenstates of the stadium billiard ($\beta=1$). Left ($k=39.045$):
 the mode has a {\em scar}
along the unstable horizontal orbit. Right ($k=39.292$): the mode is
localized in the {\em bouncing-ball region}.\label{f:scars} }
\end{center}
\end{figure}

The only counterexample to (generalized) QUE
in the case of a chaotic flow concerns the stadium billiard, or
surfaces with similar dynamical properties. This billiard admits a
1-dimensional family of marginally stable periodic orbits, the
so-called {\em bouncing-ball} orbits hitting the horizontal sides of
the stadium orthogonally (these orbits form a set of Liouville measure
zero, so they do not prevent the flow from being ergodic). In 1984
Heller \cite{Hel84} had observed that some
eigenstates are concentrated in the rectangular region (see Fig.~\ref{f:scars}).
These states were baptized {\em bouncing-ball modes}, and
studied both numerically and theoretically
\cite{HelO'Con88,BSS97}. In particular, the relative number of these modes
becomes negligible in the limit $K\to\infty$, so they are still
compatible with quantum ergodicity. Yet, no rigorous existence result
of these modes was known, until 
Hassell proved that some high-frequency eigenstates of (some)
stadia indeed have a positive mass on the bouncing ball orbits. 
To state his result, let us parametrize the shape of a
stadium billiard by the ratio $\beta$ between the length and the
height of the rectangle.
\begin{Theorem} \cite{Hassell10} For any $\eps>0$, there exists a subset $B_\eps\subset
  [1,2]$ of measure $\geq 1-4\eps$ and a number $m(\eps)>0$ such that, for any $\beta\in
  B_\eps$, the $\beta$-stadium admits a semiclassical measure with
a weight $\geq m(\eps)$ on the bouncing-ball orbits.
\end{Theorem}
Although the theorem only guarantees that a fraction $m(\eps)$ of
the semiclassical measure is localized along the bouncing-ball orbits,
the numerical studies suggest that some modes
are asymptotically fully concentrated on these orbits. Besides, such modes
are expected to exist for all ratios $\beta>0$.

\subsection{Entropy of the semiclassical measures\label{s:entropy}}

To end this section on the macroscopic properties of eigenfunctions, let us mention a
recent approach allowing to constrain the semiclassical measures occurring in a chaotic
system. This approach, initiated by Anantharaman \cite{Ana08}, consists in
proving nontrivial lower bounds for the {\em Kolmogorov-Sinai (KS) entropy} of
semiclassical measures, using the dispersion property of the
Schr\"odinger equation induced
by the hyperbolicity of the geodesic flow. The KS entropy is a common tool in classical
dynamical systems \cite{KatHas95}.
To be brief, the entropy $H_{KS}(\mu)$ of a $\Phi^t$-invariant probability measure $\mu$
is a nonnegative number which quantifies the {\em information-theoretic
complexity} of
$\mu$-typical trajectories. As such it does not direcly measure the
localization of $\mu$, but gives some information about it. Here are
some relevant properties:
\begin{itemize}
\item the measure $\mu_\gamma$ supported on a single periodic
  orbit $\gamma$  has entropy zero.
\item for an Anosov system (flow or diffeomorphism), the entropy is connected to the unstable
  Jacobian $J^u(\vx)$ (see Fig.~\ref{f:manifolds}) through the Ruelle-Pesin formula: 
\be\label{e:Ruelle}
\forall \mu\ \text{invariant,}\quad H_{KS}(\mu)\leq \int \log J^u(\vx)\,d\mu(\vx),
\ee
with equality iff $\mu$ is the Liouville measure.
\item the entropy is an affine function on the set of probability
  measures:\\
$H_{KS}\big(\alpha\mu_1+(1-\alpha)\mu_2\big) = \alpha H_{KS}(\mu_1)+(1-\alpha)H_{KS}(\mu_2)$.
\end{itemize}
In particular, the invariant measure $\alpha\mu_\gamma+(1-\alpha)\mu_L$
of a hyperbolic symplectomorphism $S$
has entropy $(1-\alpha)\lambda$, where $\lambda$ is the positive
Lyapunov exponent.

\medskip

Anantharaman considered the case of geodesic flows on manifolds $M$ of
negative curvature, see \S\ref{s:geodesic}. She proved the following
constraint on the 
semiclassical measures of $M$:
\begin{Theorem}\cite{Ana08}
Let $(M,g)$ be a smooth compact riemannian manifold of negative sectional
curvature. Then there exists $c>0$ such that any semiclassical measure
$\mu_{sc}$ of $(M,g)$ satisfies $H_{KS}(\mu_{sc})\geq c$.
\end{Theorem}
In particular, this result forbids semiclassical
measures from being convex combinations of measures $\mu_{\gamma_i}$ supported on periodic geodesics.
A more quantitative lower bound was obtained in
\cite{AnaNo07-1,AKN07}, related with the instability of the flow.
\begin{Theorem}\cite{AKN07}Under the same assumptions as above, any
  semiclassical measure must satisfy
\be\label{e:lower2}
H_{KS}(\mu_{sc})\geq \int \log J^u\,d\mu_{sc} - \frac{(d-1)\lambda_{\max}}{2},
\ee
where $d=\dim M$ and $\lambda_{\max}$ is the maximal expansion rate of
the flow.
\end{Theorem}
This lower bound was generalized to the case of the Walsh-quantized
baker's map \cite{AnaNo07-1}, and the hyperbolic
symplectomorphisms on $\t2$ \cite{Brooks10,NonCRM10}, where it takes
the form
$H_{KS}(\mu_{sc})\geq\frac{\lambda}2$. For these maps, the bound is
{\em saturated} by the half-localized semiclassical measures $\frac12(\mu_\gamma+\mu_L)$.

The bound \eqref{e:lower2} is certainly not optimal in cases of
variable curvature. Indeed, the right hand side
may become negative when the curvature varies too much. A more
natural lower bound has been obtained by Rivi\`ere in two dimensions:
\begin{Theorem}\cite{Riv08,Riv09}Let $(M,g)$ be a compact riemannian surface of
  nonpositive sectional curvature. Then any semiclassical measure
  satisfies
\be\label{e:entropy-bound}
H_{KS}(\mu_{sc})\geq \frac12\int \lambda_+\,d\mu_{sc}, 
\ee
where $\lambda_+$ is the positive Lyapunov exponent.
\end{Theorem}
The same lower bound was obtained by Gutkin for a family of
nonsymmetric baker's map \cite{Gutkin10}; he also showed that the bound is optimal for
that system. 
The lower bound \eqref{e:entropy-bound} is also expected to hold for ergodic billiards, like the
stadium; it would not contradict the existence of
semiclassical measures supported on the bouncing ball orbits.

For higher dimensional Anosov systems, one may conjecture \cite{Ana08,AnaNo07-2} the lower bound
\be\label{e:conj-bound}
H_{KS}(\mu_{sc})\geq\frac12\int \log J^u\,d\mu_{sc}\,.
\ee
Kelmer's counterexamples \cite{Kelmer06} show that
this bound may be saturated for certain Anosov
diffeomorphisms on $\IT^{2d}$. In the case of quantized
symplectomorphisms on $\IT^{2d}$ (nonnecessarily hyperbolic ones), 
Rivi\`ere \cite{Riv11} recently obtained a lower bound of the form 
$$
H_{KS}(\mu_{sc})\geq \sum_{\beta\in \Spec{S}}\max\big(\log|\beta|
-\lambda_{\max}/2 , 0\big)\,,
$$
which is sharper than \eqref{e:lower2}, but generally weaker than the
conjectured bound \eqref{e:conj-bound} in the hyperbolic (Anosov) case.

In view of the Pesin-Ruelle inequality \eqref{e:Ruelle}, we notice
that the QUE conjecture~\ref{c:QUE} amounts to improving the entropic lower
bound \eqref{e:lower2}, or even the conjecture \eqref{e:conj-bound}, to
$H_{KS}(\mu_{sc})\geq\int \log J^u\,d\mu_{sc}$.

\section{Statistical description\label{s:statistics}}

The macroscopic distribution properties described in
the previous section give a poor description of the eigenstates,
compared with our knowledge of eigenmodes of integrable systems. At
the practical level, one is interested in quantitative properties
of the eigenmodes at finite values of $\hbar$. It is also desirable to
understand their structure at the microscopic scale (the scale of the wavelength $R\sim\hbar$), or at
least some {\em mesoscopic} scale ($\hbar\ll R\ll 1$).

The results we will present are of two types. On the one hand,
{\em individual} eigenfunctions will be
analyzed statistically, e.g. by computing correlation functions or value
distributions of various representations (position density, Husimi density). 
On the other hand, one can also
perform a statistical study of a bunch of eigenfunctions (around
some large wavevector $K$), for
instance by studying how global indicators of
localization (e.g. the norms $\|\psi_j\|_{L^p}$) are distributed. 
We will not attempt to review all possible statistical
indicators, but only some ``popular'' ones.

\subsection{Chaotic eigenstates as random states?}
It has been realized quite early that the statistical properties of chaotic
eigenstates (computed numerically) could be reproduced by considering
instead
{\em random states}. The latter are, so far, the best Ansatz we
have to describe chaotic eigenstates. Yet, one should keep in
mind that this Ansatz is of a different type from the WKB Ansatz
pointwise describing individual eigenstates of integrable
systems. By definition, random states only have a chance to capture the
{\em statistical} properties of the chaotic eigenstates, but have no
chance to approximate (in $L^2$) any individual eigenfunction. This
``typicality'' of chaotic eigenstates should of course be put in
parallel with the typicality of spectral correlations, embodied by the
Random matrix conjecture  (see J.Keating's lecture).

A major open problem in quantum chaos is to prove this ``typicality''
of chaotic eigenstates. The question seems as open and difficult as the Random Matrix
conjecture.

\subsubsection{Spatial correlations}
Let us now introduce the ensembles of random states.
For simplicity we consider the Laplacian on a euclidean planar domain
$\Omega$ with chaotic geodesic flow (say, the stadium billiard). As in
 \eqref{e:Helmholtz}, we denote by $k_j^2$ the eigenvalue of $-\Delta$
corresponding to the eigenmode $\psi_j$. 

Let us recall some history.
Facing the absence of explicit expression for the
eigenstates, Voros \cite[\S 7]{Voros77} and Berry \cite{Berry77}
proposed to (brutally) approximate the Wigner measures $\mu^W_j$ of high-frequency
eigenstates $\psi_j$ by the Liouville measure $\mu_L$ on $\cE$. We've
seen that this approximation is
justified by the quantum ergodicity theorem, as long as one investigates macroscopic
properties of this Wigner function. However, Berry \cite{Berry77} also
showed that this approximation provides 
nontrivial predictions for the {\em microscopic
 correlations} of the eigenfunctions $\psi_j$. Indeed, a partial Fourier transform of
the Wigner function leads to the  {\em autocorrelation function}
describing the short-distance oscillations of $\psi$. Let us define the
correlation function by averaging over some distance $R$: in 2 dimensions,
$$
C_{\psi,R}(x,r)\defeq
\frac{1}{\pi R^2}\int_{|y-x|\leq R}dy\,\psi^*(y-r/2)\psi(y+r/2) = \frac{1}{\pi R^2} \int_{|y-x|\leq R}dy\int
dp\, e^{-irp/\hbar}\,d\mu_\psi^W(y,p)\,,
$$ 
and take $R$ to be a {\em mesoscopic scale} $k_j^{-1}\ll R \leq 1$ in
order to average over many oscillations of $\psi$. 
Inserting $\mu_L$ in the place of $\mu_\psi^W$ then provides a simple expression for this
function in the range $0\leq |r|\ll 1$:
\be\label{e:J_0}
C_{\psi,R}(x,r)\approx \frac{J_0(k|r|)}{\Vol(\Omega)}\,.
\ee
Such a homogeneous and isotropic expression could be expected from our
approximation. Replacing the Wigner distribution by $\mu_L$ suggests
that, near each point $x\in\Omega$, the
eigenstate $\psi$ is an equal mixture of particles of energy $k^2$
travelling in all possible directions.

\subsubsection{A random state Ansatz}
Yet, the Liouville measure $\mu_L$ is NOT the Wigner distribution of any quantum
state\footnote{Characterizing the function on $T^*\IR^d$ which are
  Wigner functions of individual quantum states is a nontrivial question. }. 
The next question is thus \cite{Voros77}: can one exhibit a family of quantum
states, the Wigner measures of which resemble $\mu_L$? Equivalently, the
microscopic correlations of which behave like \eqref{e:J_0}?

To account for these
isotropic correlations, Berry proposed the Ansatz of a {\em random
  superposition of isoenergetic plane waves}. One form of this Ansatz reads
\be\label{e:random}
\psi_{rand,k}(x)=\Big(\frac{2}{N\Vol(\Omega)}\Big)^{1/2}\Re 
\Big(\sum_{j=1}^N a_j \exp(k\hat n_j\cdot x)\big)\,,
\ee
where $(\hat n_j)_{j=1,\ldots,N}$ are unit vectors randomly distributed on
the unit circle, while the coefficients $(a_j)_{j=1,\ldots,N}$ are independent
identically distributed (i.i.d.) complex normal Gaussian random variables.
In order to span all possible velocity directions (within the uncertainty
principle), one should include $N\approx k$ directions $\hat n_j$
The normalization ensures that $\|\psi_{rand,k}\|_{L^2(\Omega)}\approx 1$  with high
probability when $k\gg 1$.

Alternatively, one can replace the plane waves in \eqref{e:random} by
circular-symmetric waves, namely Bessel functions. The random state
then reads
\be\label{e:random2}
\psi_{rand,k}(r,\theta)=(\Vol(\Omega))^{-1/2}\sum_{m=-M}^M b_m J_{|m|}(kr)\,e^{i m \theta}\,,
\ee
where the coefficients $b_m$ are i.i.d. complex Gaussian satisfying the symmetry
$b_m=b_{-m}^*$, and $M\approx k$.
Both these random ensembles asymptotically produce the same
statistical results when $k\to\infty$.

The random state $\psi_{rand,k}$ satisfies the equation $(\Delta+k^2)\psi=0$ in the
interior of $\Omega$. Furthermore, it
satisfies a ``local quantum ergodicity'' property: for any observable
$f(x,p)$ supported in the interior of $T^*\Omega$, 
the quantum average
$\la\psi_{rand,k},\hat f_{k^{-1}}\psi_{rand,k}\ra\approx \mu_L(f)$ with
high probability (more is known about the statistics of these arerages, see
\S\ref{s:distrib}).

The stronger claim is that, in the interior of $\Omega$, the  {\em
  local statistical properties} of $\psi_{rand,k}$, including
microscopic ones,
should be similar with those of the eigenstates $\psi_j$ with
wavevectors $k_j\approx k$.  

The correlation functions of eigenstates of chaotic planar billiards
have been numerically studied, and compared with the random models, 
see e.g. \cite{McDoKauf88,BackSchub02}. The agreement with \eqref{e:J_0} is
fair for some eigenmodes, but not so good for others; 
in particular the authors observed some anisotropy in
the ``true'' correlation functions, which may be related to some form of {\em scarring}
(see \S\ref{s:scars}), or to the bouncing-ball modes of the stadium billiard.

The {\em value distribution} of the random wavefunction \eqref{e:random} is Gaussian, and
compares very well with numerical studies of eigenmodes of chaotic billiards
\cite{McDoKauf88}. A similar analysis has been performed for
eigenstates of the Laplacian on a compact surface of constant negative curvature
\cite{AurStein93}. In this geometry the random Ansatz was defined in
terms of adapted circular hyperbolic waves. The authors checked that the coefficients of the
individual eigenfunctions, when expanded in these hyperbolic waves, are indeed Gaussian
distributed. They also checked that the value distribution of 
individual eigenfunctions $\psi_j(x)$ is Gaussian to a good accuracy, without any
exceptions.

\subsubsection{On the distribution of quantum averages \label{s:distrib}}
The random state model
also predicts the statistical distribution of the quantum averages
$\la \psi_j,\hat f_\hbar\psi_j\ra$, namely the average of the observable $f$ w.r.to the
Wigner distributions, $\mu_j^W(f)$. The quantum variance estimate in
the proof of Thm.~\ref{thm:QE} shows
that the distribution of these averages becomes semiclassically concentrated around the
classical value $\mu_L(f)$. Using a mixture of
semiclassical and random matrix theory arguments, Feingold and Peres
\cite{FeinPer86} conjectured that, in the semiclassical 
limit, the quantum averages of eigenstates in a small energy window
$\hbar k_j\in [1-\eps,1+\eps]$ should be Gaussian distributed, with
mean $\mu_L(f)$ and variance related with the {\em
  classical variance} of $f$. The latter is
defined as the integral of the autocorrelation function $C_{f,f}(t)$
(see \eqref{e:mixing}):
$$
\Var_{cl}(f)= \int_{\IR} C_{f,f}(t)\,dt\,.
$$
A more precise semiclassical derivation \cite{Ekh+95}, using the Gutzwiller trace formula, and
supported by numerical computations on several chaotic systems, 
confirmed both the Gaussian distribution of the quantum averages, and
showed the following connection between quantum and classical
variances (expressed in semiclassical notations):
\be
\Var_{\hbar}(f)\sim g\,\frac{\Var_{cl}(f)}{T_H}\,.
\ee
Here $g$ is a symmetry factor
($g=2$ in presence of time reversal symmetry, $g=1$ otherwise), and
$T_H=2\pi\hbar\bar\rho$ is the Heisenberg time, where $\bar\rho$ is the
smoothed density of states. In the case of the semiclassical
Laplacian on a compact surface or a planar domain,
the mean density $\bar{\rho}$ is drawn from Weyl's law \eqref{e:Weyl}, so we get
$$
\Var_{\hbar}(f)\sim g\hbar\Var_{cl}(f)/\Vol(\Omega)\,
$$ 
Equivalently, the quantum variance corresponding to wavevectors
$k_j\in[K,K+1]$ is predicted to take the value (when $K\gg 1$):
\be\label{e:qu-var-decay}
\Var_{K}(f)\sim \frac{g}{K}\,\frac{\Var_{cl}(f)}{\Vol(\Omega)}\,.
\ee
Successive numerical studies on chaotic euclidean billiards
\cite{BSS98,Bar06} and manifolds or billiards of negative curvature
\cite{AurTag98} globally
confirmed this prediction for the quantum variance, as well as the
Gaussian distribution of the quantum averages at high
frequency. Still, the
convergence to this law can be slowed down for billiards admitting bouncing-ball
eigenmodes, like the stadium billiard \cite{BSS98}.

For a generic chaotic system, rigorous semiclassical methods could only prove 
logarithmic upper bounds for the quantum variance \cite{Schu06},
$$
\Var_{\hbar}(f)\leq C/|\log\hbar|\,.
$$
Schubert showed that this slow decay can be sharp
for certain eigenbases of the quantum cat map, in the case of large
spectral degeneracies \cite{Schu08} (as we have seen in
\S\ref{s:counter}, such degeneracies are also
responsible for the existence exceptionally localized eigenstates, so
a large variance is not surprising).

The only systems for which an algebraic decay is proved
are of arithmetic nature. Luo and Sarnak \cite{LuoSar04}, and then
Zhao \cite{Zhao10}
proved that, in the case of the modular domain $M=SL_2(\IZ)\backslash
\IH$ (a noncompact, finite volume arithmetic surface for which the
Laplacian admits many $L^2$ eigenstates), the quantum
variance corresponding to high-frequency Hecke eigenfunctions is of the form 
$\Var_{K}(f)=\frac{B(f)}{K}$: the polynomial decay is the same as in
\eqref{e:qu-var-decay}, but the coefficient $B(f)$ is equal
the classical variance ``decorated'' by an extra factor of arithmetic
nature.

More precise results were obtained for
quantum symplectomorphisms on the 2-torus. Kurlberg and Rudnick \cite{KurRud05} studied
the distribution of quantum averages
$\{\sqrt{N}\la \psi_{N,j},\hat f_N\psi_{N,j}\ra,\ j=1,\ldots,N\}$, where the $(\psi_{N,j})$ form a Hecke eigenbasis
of $U_N(S)$ (see \S\ref{s:arithmetic}). In the semiclassical limit
$N\gg 1$, the
variance is asymptotically of the form
$\frac{B(f)}{N}$, with $B(f)$ an ``arithmetically decorated'' classical variance. They
also computd the fourth moment of the distribution, which suggests that the latter is not Gaussian, but given
by a combination of several semicircle laws on $[-2,2]$ (or Sato-Tate
distributions). The fifth moment was recently computed by Rosenzweig,
and shown compatible with this conjecture \cite{Rosen10}.
The same semicircle law had been shown in \cite{KurRud01b} to correspond to
the asymptotic distribution of the (position) coefficients of the Hecke eigenstates, at least
for $N$ along a subsequence of ``split primes''.

\subsubsection{Maxima of eigenfunctions}
Another interesting quantity is the statistics of the maximal
values of eigenfunctions, that is their $L^\infty$ norms, or more generally their $L^p$
norms for $p\in(2,\infty]$  (we always assume
the eigenfunctions to be $L^2$-normalized). The maxima belong to the far
tail of the value distribution, so their behaviour is a priori
uncorrelated with the Gaussian nature of the latter. 

The random wave model gives the following estimate
\cite{AuBaSchuTag99}: for $C>0$ large enough,
\be\label{e:sup-bound}
\frac{\|\psi_{rand,k}\|_{\infty}}{\|\psi_{rand,k}\|_{2}}\leq
C\sqrt{\log k}\quad \text{with high probability when $k\to\infty$}
\ee
Numerical tests on some euclidean chaotic billiards and a
surface of negative curvature show that this order of magnitude is
correct for chaotic eigenstates \cite{AuBaSchuTag99}.
Small variations were observed between arithmetic/non-arithmetic
surfaces of constant negative curvature, the sup-norms appearing
slightly larger in the arithmetic case, but still compatible with
\eqref{e:sup-bound}. For the euclidean billiards, the largest maxima
occured for states {\em scarred} along a periodic orbit (see \S\ref{s:scars}).

Mathematical results concerning the maxima of eigenstates of
generic manifolds of negative curvature are scarce. A general upper bound
\be\label{e:Hormander}
\|\psi_j\|_\infty\leq C\,k_j^{(d-1)/2}
\ee
holds for arbitrary compact manifolds \cite{Horm68}, and is saturated in the
case of the standard spheres.
On a manifold of negative curvature, this upper bound can be improved
by a factor $(\log k_j)^{-1}$, taking into account a better
bound on the remainder in Weyl's law.

Once again, more precise results have been obtained only in the case of
Hecke eigenstates on arithmetic manifolds.
Iwaniec and Sarnak \cite{IwaSar95} showed that,
for a certain family of arithmetic surfaces, the Hecke eigenstates
satisfy the bound 
$$
\|\psi_j\|_\infty\leq C_\eps\,k_j^{5/12+\eps}\,,
$$
and conjecture a bound 
\be\label{e:conj-sup}
\|\psi_j\|_\infty\leq C_\eps\,k_j^{\eps}\,,
\ee
slightly larger than for the 
random wave model. More recently, Mili\'cevi\'c \cite{Mili09} showed that,
on certain arithmetic surfaces, a subsequence
of Hecke eigenstates satisfies a {\em lower} bound
$$
\|\psi_j\|_\infty\geq C\exp\Big\{\left(\frac{\log k_j}{\log\log k_j}\right)^{1/2}(1+o(1))\Big\}\,,
$$
thereby violating the random wave model, but still compatible with the
conjecture \eqref{e:conj-sup}. The large values are
reached on specific CM-points of the surface, of arithmetic nature.

On certain higher dimensional arithmetic manifolds, Rudnick and Sarnak \cite{RudSar94} had already
identified some Hecke eigenstates with larger values, namely
$$
\|\psi_j\|_\infty\geq C\,k_j^{1/2}\,.
$$
A general discussion of this phenomenon appears in the recent
work of Mili\'cevi\'c \cite{Mili10}; the author presents a larger family
of arithmetic 3-manifolds featuring eigenstates with abnormally large
values, and conjectures that his list is exhaustive.
\begin{figure}[htbp]
\begin{center}
\includegraphics[width=.95\textwidth]{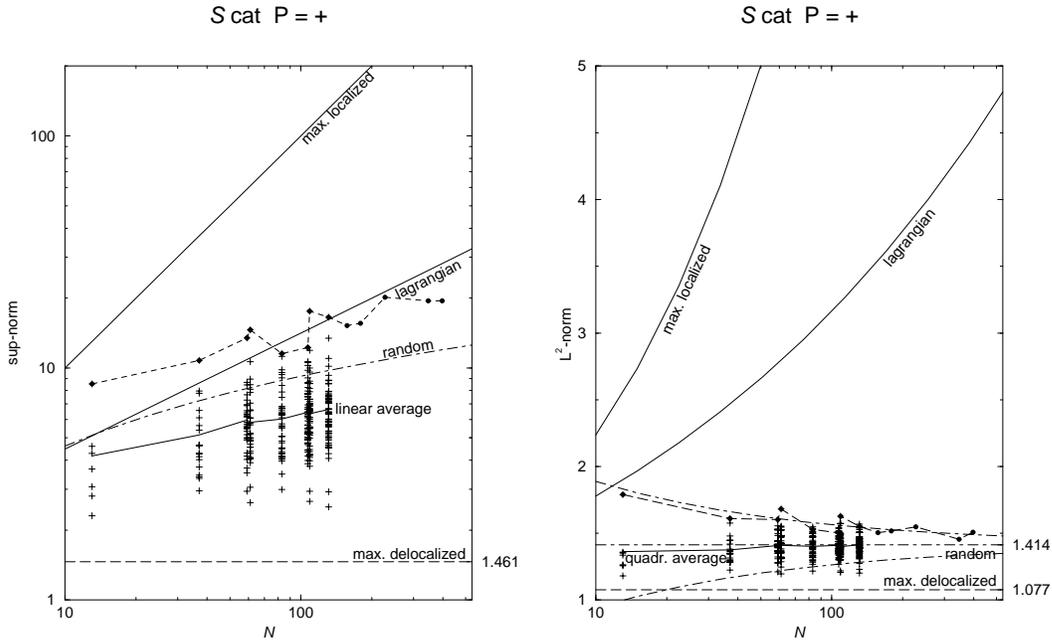}
\caption{$L^\infty$ (left) and $L^2$ (right) norms of the Husimi densities for
  eigenstates of the quantum symplectomorphism $U_N(S_{DEGI})$ (crosses; the
  dots indicate the states maximally scarred at the origin). The
  data are compared with the values for maximally
  localized states (Fig.~\ref{f:integ}, right), lagrangian states
  (Fig.~\ref{f:integ}, center), 
random states (Fig.~\ref{f:cat-Hecke}, bottom right) and
  the maximally delocalized state. (Repr. from \cite{NonVor98}). \label{f:cat-Husimi-stats}}
\end{center}
\end{figure}

\subsubsection{Random states on the torus\label{s:random-t2}}

In the case of quantized chaotic maps on $\t2$, one can easily setup
a model of random states mimicking the statistics of eigenstates.
The choice is particularly simple when the map does not possess any particular
symmetry: the ensemble of
random states in $\hn$ is then given by
\be\label{e:random-t2}
\psi_{rand,N}=\frac{1}{\sqrt{N}}\sum_{\ell=1}^N a_\ell\,e_\ell\,,
\ee 
where $(e_\ell)_{\ell=1,\ldots,N}$ is the orthonormal basis \eqref{e:e_j} of $\hn$, and the
$(a_\ell)$ are i.i.d. normal complex Gaussian variables. This random
ensemble is $U(N)$-invariant, so it takes the same form w.r.to any
orthonormal basis of $\hn$. For this ensemble, the Husimi
function $\cH_{\psi_{rand,N}}(\vx)$ admits, at each point $\vx\in \t2$, an
exponential value distribution with average $1+f_N(\vx)$, where $f_N$ is
  exponentially small when $N\to\infty$.

This random model, and some variants taking into account symmetries, 
have been used to describe the spatial, but also the phase space
distributions of eigenstates of quantized chaotic maps. In
\cite{NonVor98}, various indicators of localization of the
Husimi densities \eqref{e:Husimi}\footnote{By construction the Husimi
  densities are $L^1$-normalized.} have been computed for this random model, and compared
with numerical results for the eigenstates of the quantized ``cat''
and baker's maps (see
Fig.~\ref{f:cat-Husimi-stats}). 
The distributions seem compatible with the random state model, except
for the large deviations of the sup-norms of the Husimi
densities, due to eigenstates ``scarred'' at the fixed point at the origin
(see Fig.~\ref{f:cat-Hecke}). 

In  \cite{KurRud01b} the Hecke eigenstates of $U_N(S)$, expressed in the position basis
$(e_\ell)$ as in \eqref{e:random-t2}, were shown to satisfy nontrivial $\ell^\infty$
bounds\footnote{In the chosen normalization, the trivial bound reads
  $\|\psi\|_\infty\leq N^{1/2}$.}:
$$
\|\psi_{N,j}\|_\infty= \sup_{\ell}|a_\ell|\leq C_\eps \,N^{3/8+\eps}.
$$
For $N$ along a subsequence of ``split primes'', the
description of the individual
Hecke eigenstates is much more precise. Their position coefficients are
uniformly bounded, $\|\psi_{N,j}\|_\infty\leq 2$, and the value
distribution of individual eigenstates $\{|\psi_{N,j}(\ell/N)|,\
\ell=0,\ldots,N-1|\}$ is asymptotically given by the
semicircle law on $[0,2]$, showing that these eigenstates are
very different from Gaussian random states. 
In spite of this fact, the value distribution of the Husimi function
$\cH_{\psi_{N,j}}(\vx)$, which, at each $\vx$ involves $\approx \sqrt{N}$
position coefficients, appears to be exponential, like in the random
state model.

\subsection{Scars of periodic orbits \label{s:scars}}

\begin{figure}[htbp]
\begin{center}
\includegraphics[angle=-90,width=.9\textwidth]{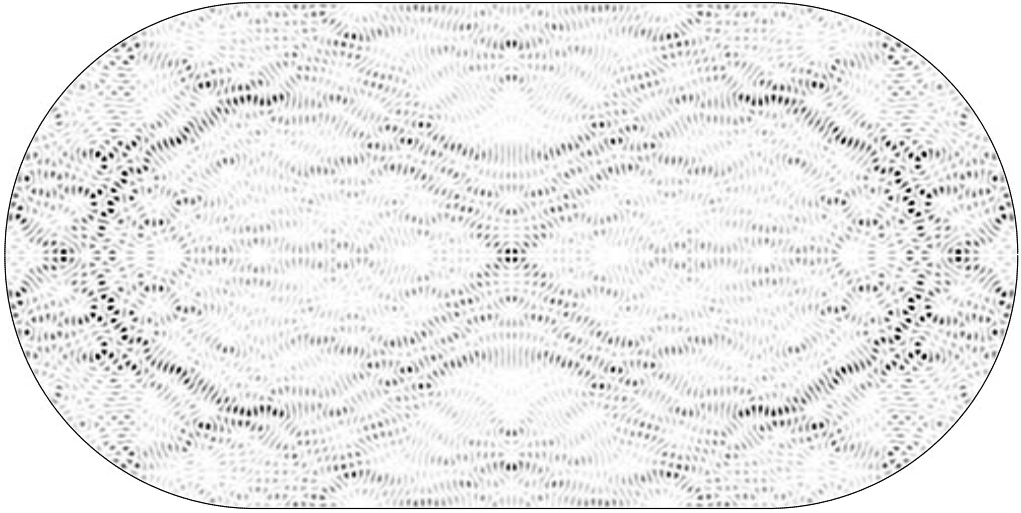}
\end{center}
\caption{A high-energy eigenmode of the stadium billiard ($k\approx
  130$). Do you see any scar?}\label{f:eigen130}
\end{figure}

Around the time the quantum ergodicity theorem was proved, an
interesting phenomenon was observed by Heller in
numerical studies on the stadium billiard \cite{Hel84}. He
noticed that for certain eigenfunctions, the spatial density
$|\psi_j(x)|^2$ is {\em abnormally enhanced}
along one or several {\em unstable} periodic geodesics. He called such an enhancement a {\em scar} of
the periodic orbit on the eigenstate $\psi_j$. See
Fig.~\ref{f:scars} (left) and Fig.~\ref{f:eigen130} for scars at low and relatively high frequencies.

This phenomenon was observed to persist at higher and higher frequencies for
various euclidean billiards \cite{VS95,AuBaSchuTag99}, but was
not detected on manifolds of negative curvature \cite{AurStein93}.
Could scarred states represent counterexamples to quantum unique ergodicity?
More precise numerics \cite{Bar06} showed that
the weights of these enhancements near the relevant periodic
orbit decay in the
high-frequency limit, because the {\em areas} covered by these enhancements decay
faster than their {\em intensities}. As a result, a sequence of scarred states may still
become equidistributed in the classical limit \cite{KapHel98,Bar06}.

In order to quantitatively characterize the scarring
phenomenon, it appeared more convenient to switch to 
phase space representations, in particular the Husimi
density of the boundary function $\Psi(q)=\partial_\nu\psi(q)$, which lives on
the phase space $T^*\partial\Omega$ of the billiard map
\cite{CresPerCha93,TuaVor95}. In this representation, scars were detected as
enhancements of the Husimi density $\cH_{\Psi_j}(\vx)$ on periodic
phase space points $\vx\in T^*\partial\Omega$ of the billiard map (see
Fig.~\ref{f:Poincare}, bottom). 

Similar studies were performed in the case of quantum chaotic maps on
the torus, like the baker's map
\cite{Sar90} or hyperbolic symplectomorphisms \cite{NonVor98}. The 
scarred state showed
in Fig.~\ref{f:cat-Hecke} (left) has the largest value of the Husimi
density among all eigenstates of $U_N(S_{DEGI})$, but it is nevertheless a Hecke
eigenstate, so that its Husimi measure should be (macroscopically)
close to $\mu_L$. In this example, the scarring phenomenon is a
microscopic (or mesoscopic) phenomenon, compatible
with quantum unique ergodicity. This sequence of scarred eigenstates
also featured abnormally large values of the Husimi density (see
Fig.~\ref{f:cat-Husimi-stats}, left).

\subsubsection{A statistical theory of scars}

Heller first tried to explain the scarring phenomenon using 
the {\em smoothed local density of states}
\be\label{e:LDOS}
S_{\chi,\vx_0}(E)=\sum_j \chi(E-E_j)|\la
\varphi_{\vx_0},\psi_j\ra|^2=\la \varphi_{\vx_0},\chi(E-\hat H_\hbar)\varphi_{\vx_0}\ra\,,
\ee
where $\varphi_{\vx_0}$ is a Gaussian wavepacket \eqref{e:wavepacket} sitting on a
point of the periodic orbit, and $\chi\in C^\infty_c(\IR)$ an energy cutoff.
The function $\chi$ is constrained by the fact that this expression is
estimated through its Fourier transform, namely
the time autocorrelation function
\be\label{e:autocorr}
t\mapsto \la \varphi_{\vx_0},U_\hbar^t\varphi_{\vx_0}\ra\,\tilde\chi(t)\,,
\ee
where $\tilde\chi$ is the $\hbar$-Fourier transform of $\chi$. Because
we can control the evolution of $\varphi_{\vx_0}$ only up to the
Ehrenfest time \eqref{e:Ehrenfest}, we must take $\tilde\chi$
supported on the interval $[-T_E/2,T_E/2]$, so that $\chi$ is
essentially supported on an interval of
width $\gtrsim \frac{\hbar}{|\log\hbar|}$.
Since $U_\hbar^t\,\varphi_{\vx_0}$ comes back to the point $\vx_0$ at
each period $T$, $S_{\chi,\vx}(E)$ has peaks at the {\em Bohr-Sommerfeld energies} of
the orbit, separated by $2\pi\hbar/T$ from one another. However, due to the hyperbolic
spreading of the wavepacket, these peaks have widths
$\sim\lambda\hbar/T$, where $\lambda$ is the Lyapunov exponent of the
orbit. Hence, the peaks can only be significant for small enough $\lambda$,
that is weakly unstable orbits.
Even if $\lambda$ is small, the width $\lambda\hbar/T$ becomes much larger than the mean level
spacing $1/\bar\rho\sim C\hbar^d$ in the semiclassical
limit, so that each peak encompasses
many eigenvalues. In particular, this mechanism 
can not predict {\em which} individual eigenstate will show an
enhancement at $\vx_0$, nor can it predict the value of the enhancements.

Following Heller's work, Bogomolny \cite{Bogo88} and Berry \cite{Berry89} showed that
certain linear combinations of nearby eigenstates have an
``extra intensity'' in the spatial density (resp. ``oscillatory
corrections'' in the Wigner density) around a certain
number of closed geodesics. In the
semiclassical limit, these combinations also involve many eigenstates
in some energy window.

\medskip

A decade later, Heller and Kaplan developed a ``nonlinear'' theory of scarring,
which proposes a {\em statistical} definition of this phenomenon \cite{KapHel98}. They noticed that, given an
energy interval $I$ of width $\hbar^2\ll|I|\ll\hbar$, the distribution of the overlaps
$\{|\la\varphi_{\vx},\psi_j\ra|^2,\,E_j\in I\}$ depends on the phase point
$\vx$: if $\vx$ lies on a (mildly unstable) periodic orbit, the distribution is
spread between some large values (scarred states) and some low values
(antiscarred states). On the opposite, if $\vx$ is a ``generic''
point, the distribution of the overlaps is narrower. 

This remark was made quantitative by defining a {\em stochastic model} for the unsmoothed
local density of states $S_\vx(E)$ (that is, taking $\chi$ in
\eqref{e:LDOS} to be a delta function), as an effective way to take into account
the (uncontrolled) long time recurrences in the autocorrelation
function \eqref{e:autocorr}.
According to this model, the overlaps
$|\varphi_{\vx},\psi_j\ra|^2$ in an energy window $I$ should behave
like random exponential variables, of {\em expectation given by the smoothed
local density} $S_{\chi,\vx}(E)$. Hence, if $\vx$ lies on a short
periodic orbit, the eigenstates of energies
close to the Bohr-Sommerfeld energies (where $S_{\chi,\vx}(E)$ is maximal)
statistically have larger overlaps with $\varphi_{\vx}$, while states
with energies close to the anti-Bohr-Sommerfeld energies
statistically have smaller overlaps. The concatenation of these
exponential random variables with smoothly varying expectations produces a
non-exponential distribution, with a tail larger than the one predicted
by Berry's random model.
On the opposite, if $\vx$ is a ``generic'' point, the expectation should
not depend on the energy, and the full distribution of the $|\varphi_{\vx},\psi_j\ra|^2$
remains exponential.

Although not rigorously justified, this statistical definition of scarring
gives quantitative predictions, and can be viewed as an interesting ``dynamical
correction'' of the random state model \eqref{e:random}.

\section{Nodal structures\label{s:nodal}}

For a moment let us focus on the case of a Laplacian $\Delta_\Omega$ on a
euclidean billiard, or a compact riemannian manifold. 
After having described the ``macroscopic skeleton'' of the
eigenfunctions $\psi_j$, namely their semiclassical measures, and the possible large
values taken by $|\psi_j(x)|$ near a periodic orbit or elsewhere, we
now focus on the opposite feature of these eigenfunctions, namely their {\em nodal
sets} $\cN_{\psi_j}=\{x\in\Omega,\,\psi_j(x)=0\}$. Since the
eigenfunctions $\psi_j$ can be chosen real, their nodal set is a union
of hypersurfaces (in 2 dimensions, nodal lines), which sometimes
intersect each other, or intersect the boundary $\partial\Omega$. 
This set separates the connected domains
where $\psi_j(x)$ has a definite sign, called {\em nodal domains}. The
nodal set can be viewed as a {\em microscopic skeleton} of the
eigenfunction $\psi_j$: it fully determines the function (up to
a global factor), and the typical
scale separating two nearby hypersurfaces is the wavelength
$k_j^{-1}$ (or $\hbar_j$ in the semiclassical formalism). 

The study of the nodal patterns of eigenfunctions has a long history in
mathematical physics and riemannian geometry. Except for integrable
systems\footnote{Actually, nodal domains are well-understood only for
  {\em separable} systems, a stronger assumption than integrability.},
we have no explicit knowledge of these sets. However, some global
properties are known, independently on any assumption on the
geometry. In 1923 Courant \cite{Courant23} showed that, for
the Dirichlet Laplacian on a planar
domain $\Omega\subset \IR^2$, the number of nodal domains $\nu_j$ of the
$j$-th eigenstate (counted with multiplicities) satisfies $\nu_j\leq
j$, an inequality which is an equality in
dimension 1. This upper bound was
sharpened by Pleijel \cite{Pleijel56} for high-frequency eigenstates:
$$
\limsup_{j\to\infty}\frac{\nu_j}{j}\leq 0.692\,.
$$

\subsection{Nodal count statistics for chaotic eigenstates}
\begin{figure}
\begin{center}
\includegraphics[width=0.4\textwidth]{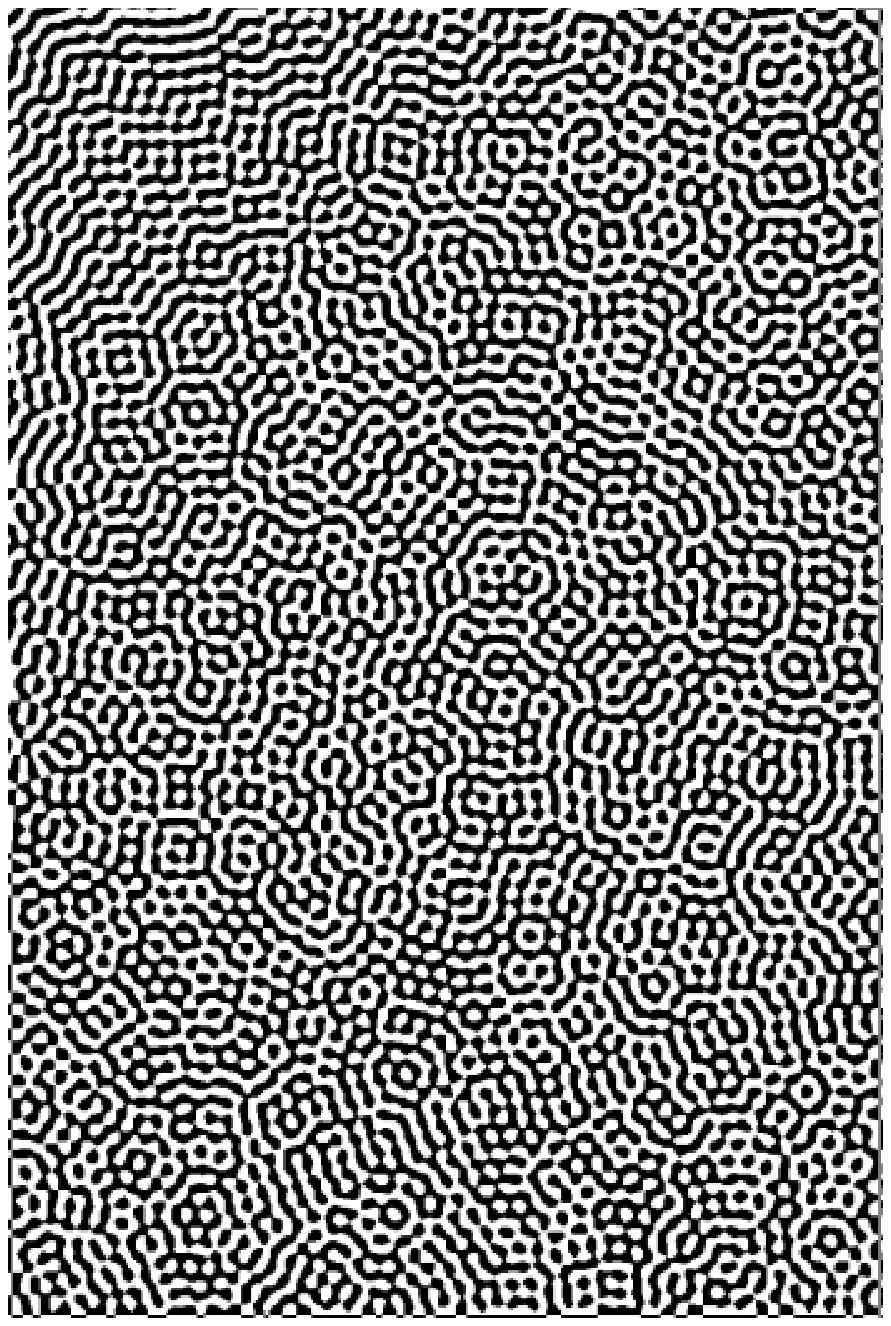}\hspace{1cm}
\includegraphics[width=0.5\textwidth]{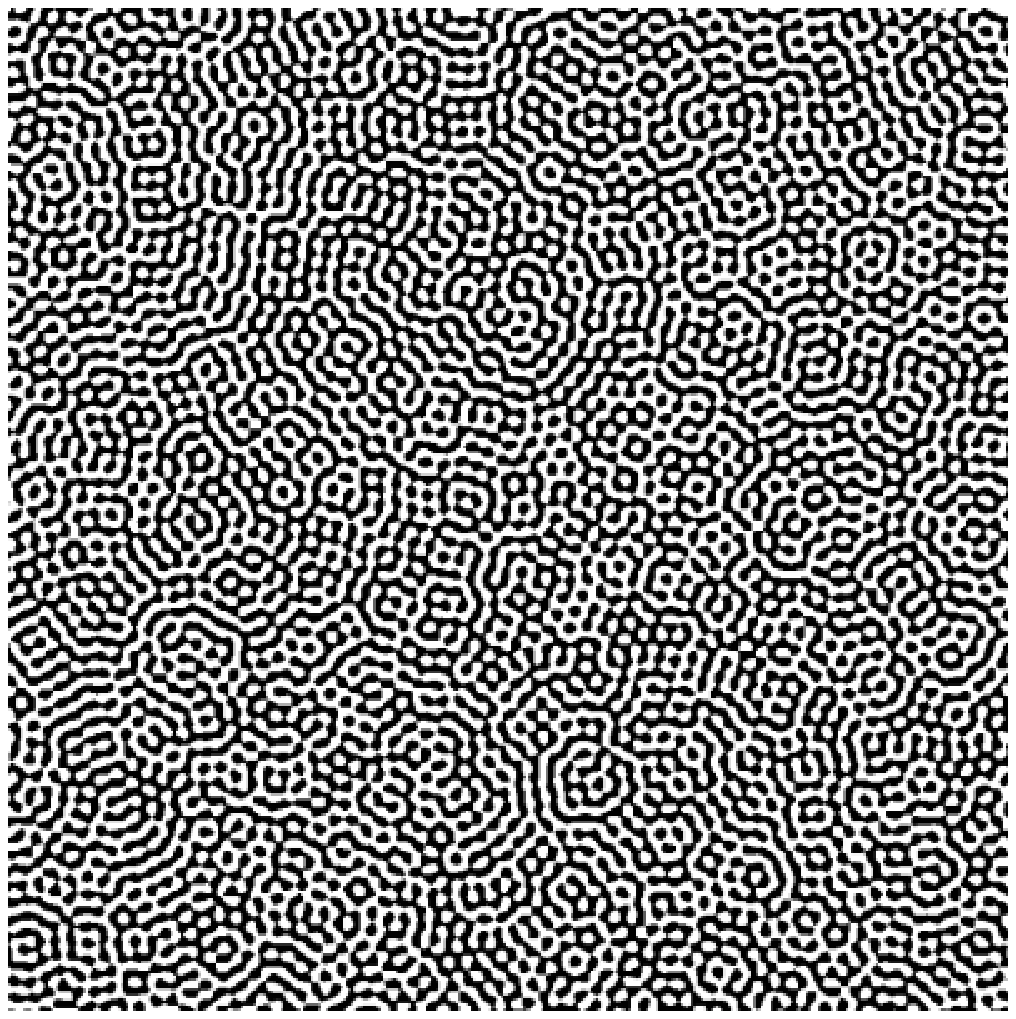}
\caption{Left: nodal domains of an eigenfunction of the
  quarter-stadium with $k\approx 100.5$ (do you see the boundary of
  the stadium?). Right: nodal domains of a random state
  \eqref{e:random2} with $k=100$ (Repr. from \cite{BogoSchm07}).}\label{f:nodal}
\end{center}
\end{figure}
The specific study of nodal structures of eigenstates for {\em
  chaotic} billiards is more recent.
Blum, Gnutzmann and Smilansky \cite{BlumGnutSmil02} seem to be the
first authors using nodal statistics to differentiate regular
from chaotic wavefunctions.
They compared
the {\em nodal count sequence} $(\xi_j=\frac{\nu_j}{j})_{j\geq 1}$ for separable vs.
chaotic planar domains, and observed different statistical behaviours. In
the separable case the nodal lines can be explicitly computed: they
form a ``grid'', defined in terms of the two independent invariant actions
$I_1,I_2$. The distribution of the sequence is peaked near some
value $\xi_m$ depending on the geometry. 

In the chaotic case, the authors found that very few nodal lines
intersect each other, and conjectured that the sequence 
should have the same statistics as in the case of the random
models (\ref{e:random},\ref{e:random2}), therefore showing some {\em
  universality}. The numerical plot of Fig.~\ref{f:nodal} perfectly
illustrates this assertion. In the random model, the
random variable corresponding to $\xi_j=\frac{\nu_j}{j}$ should be
\be
\frac{\nu(k)}{\bar N(k)},\qquad k\gg 1\,,
\ee
where $\nu(k)$ is the number of nodal
domains in $\Omega$ for the random state \eqref{e:random}, and $\bar
N(k)=\frac{\Vol(\Omega)k^2}{4\pi}$ is the integrated density
of states in $\Omega$.

\begin{figure}
\begin{center}
\includegraphics[width=0.4\textwidth]{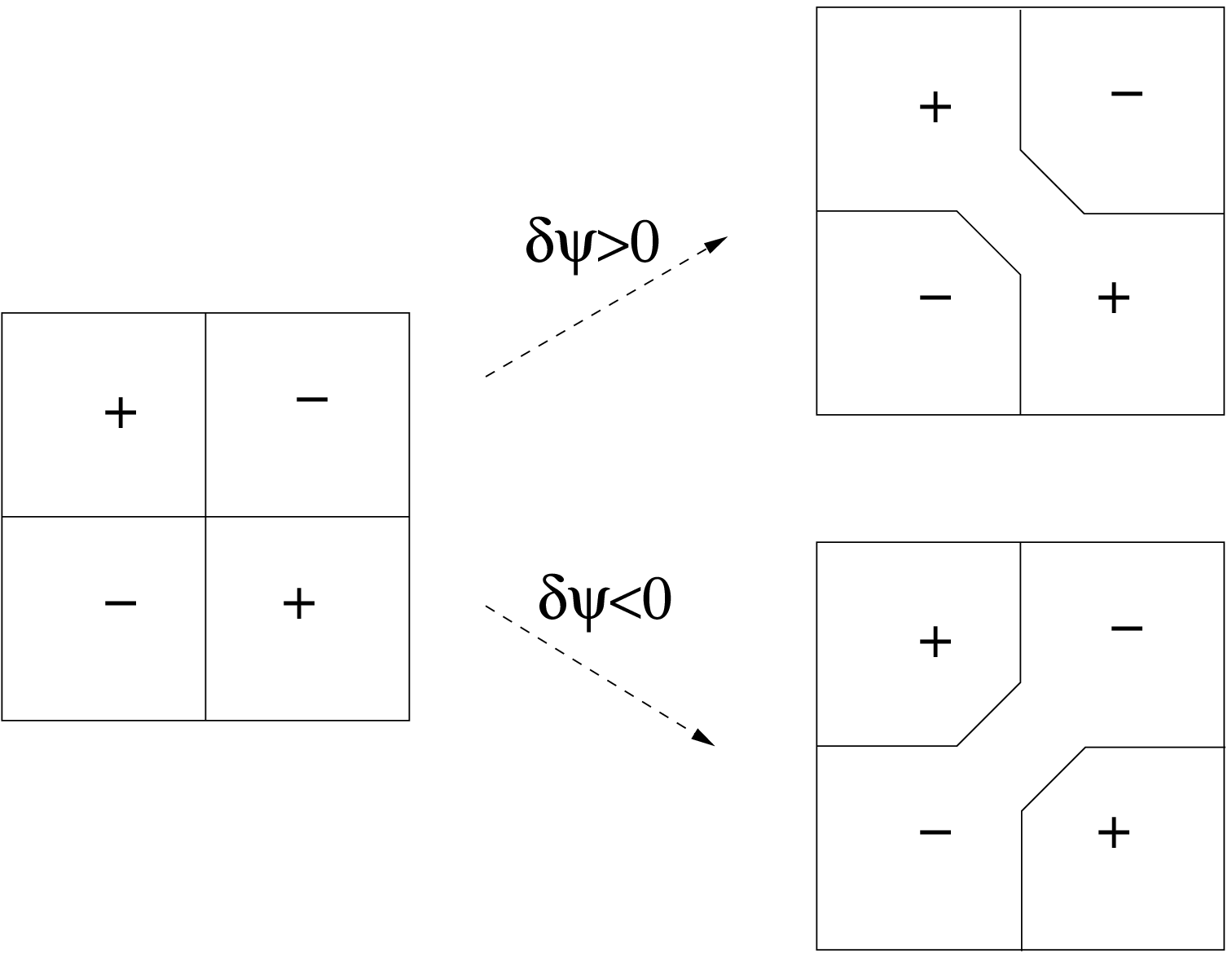}\hspace{.5cm}
\includegraphics[width=0.5\textwidth]{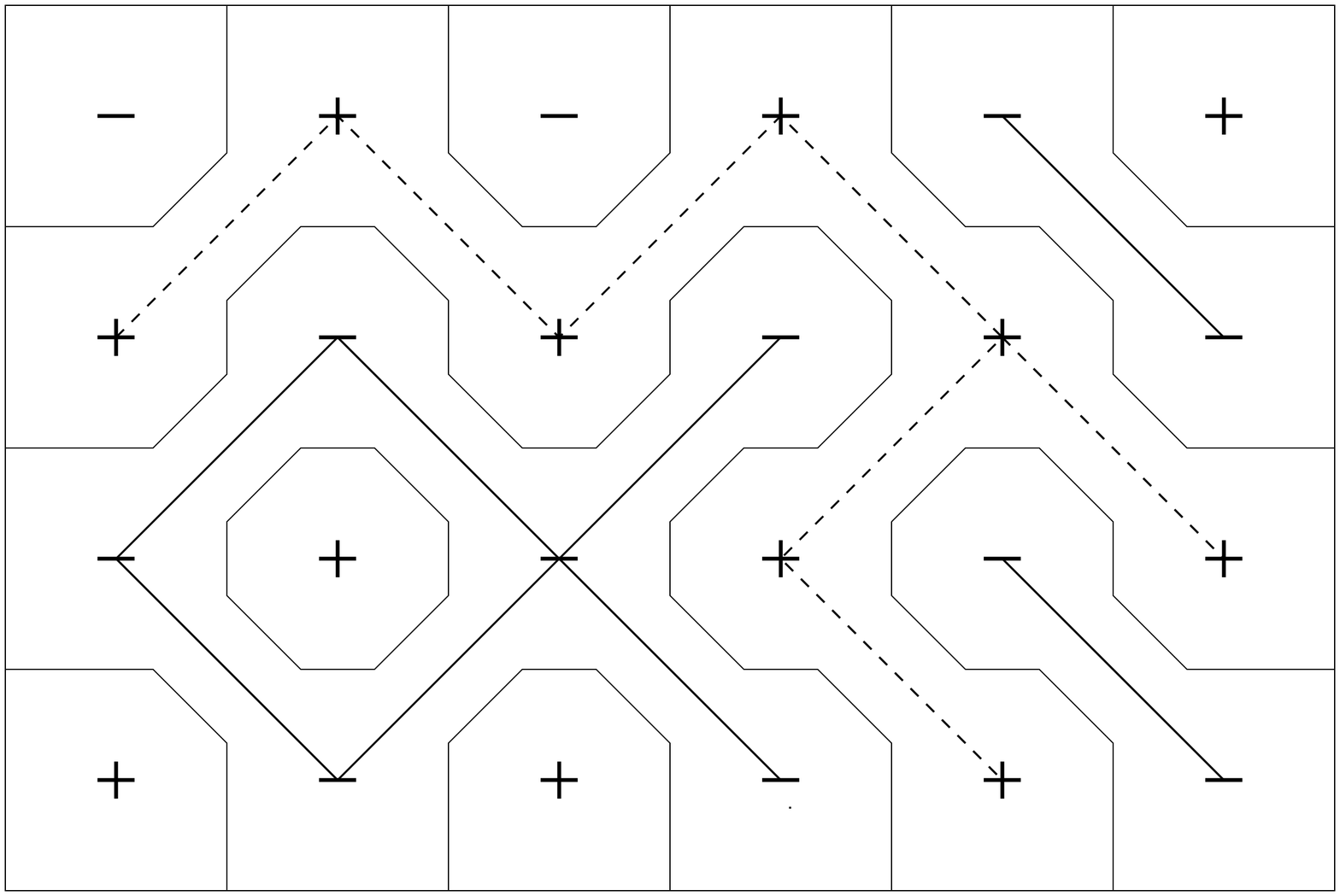}
\caption{Construction of the random-bond percolation model. Left: starting
  from a square grid where $\psi(x)$ alternatively takes positive $(+)$ and
  negative $(-)$ values, a small perturbation $\delta\psi(x,y)$ near
  each crossing creates an avoided crossing. Right: the resulting
  positive and negative nodal domains can be described by drawing
  bonds (thick/dashed lines) between adjacent sites (reprinted from
  \cite{BogoSchm02}). \label{f:perco}}
\end{center}
\end{figure}

Motivated by these observations, Bogomolny and Schmit proposed a heuristic
{\em percolation model} to compute the nodal count statistics in random or chaotic
wavefunctions \cite{BogoSchm02}. 
Their model starts from a separable eigenfunction, of the form 
$\psi_0(x,y)=\cos(kx/\sqrt{2})\cos(ky/\sqrt{2})$,
for which nodal lines form a grid; they perturb this function near
each intersection, so that each crossing becomes an avoided crossing
(see Fig.~\ref{f:perco}). Although the length of
the nodal set is almost unchanged, the structure of the {\em nodal domains}
is drastically modified by this perturbation.
Assuming these local perturbations are {\em uncorrelated}, they obtain
a representation of the nodal domains as {\em clusters} of a critical bond percolation
model, a well-known model in 2-dimensional statistical mechanics. 

The (somewhat amazing) claim made in \cite{BogoSchm02} is that such a
perturbation of the separable wavefunction has the same nodal count
statistics as a random function, eventhough the latter is
very different from the former in many respects (e.g their value distribution). 
The underlying idea is the following: the uncorrelated local
perturbations $\delta\psi(x,y)$ instantaneously transform the
microscopic square nodal domains into mesoscopic
(sometimes macroscopic) ``fractal'' domains, which then do not vary
too much upon strengthening the perturbation such as to reach a
typical random state.
The high-frequency limit ($k\to\infty $) for the random state
(\ref{e:random},\ref{e:random2}) corresponds to the
thermodynamic limit in the percolation model, a limit in which
the statistical properties of percolation clusters have been much
investigated. 
The nodal count ratio \eqref{e:nodal-count-random} counts the number of
clusters on a lattice of $N_{tot}=\frac{2}{\pi}\bar N(k)$ sites.
The distribution of the
number of domains/clusters $\nu(k)$ was computed in this limit \cite{BogoSchm02}:
it is a Gaussian with asymptotic properties
\be\label{e:nodal-count-random}
\frac{\la \nu(k)\ra}{\bar N(k)}\xrightarrow{k\to\infty} 0.0624,\qquad 
\frac{\Var(\nu(k))}{\bar N(k)}\xrightarrow{k\to\infty} 0.0502\,.
\ee
Nazarov and Sodin \cite{NazSod09} have considered random spherical harmonics on
the 2-sphere (namely, Gaussian random states within each $2n+1$-dimensional
eigenspace), and proved that when $n\to\infty$ the distribution of the
number of nodal domains of $\psi_{rand,n}$ becomes peaked near the value $a\,n^2$,
for some constant $a>0$. Although they 
compute neither the constant $a$, nor the variance of the distribution, this result 
indicates that the percolation model may indeed correctly predict the nodal
count statistics for Gaussian random states.

\begin{rem}We now have {\em two} levels of modelization.
First, the chaotic eigenstates are statistically
modelled by the random states
(\ref{e:random},\ref{e:random2}).
Second, the nodal structure of random states is modelled by 
critical percolation. These two conjectures appeal to different
methods: the second one is a purely statistical problem, while the first one
belongs to the ``chaotic=random'' meta-conjecture.
\end{rem}

\subsection{Other nodal observables}

According to the Bogomolny-Schmit percolation mode, the statistical distribution of the {\em
  areas} of nodal domains has 
the form $\cP(\cA)\sim \cA^{-187/91}$. 
Of course, this scaling can only hold in the mesoscopic range $k^{-2}\ll \cA\ll 1$,
  since any domain has an area $\geq C/k^2$. 

The number $\tilde\nu_j$ of nodal lines of $\psi_j$ {\em touching the
  boundary} $\partial\Omega$ is also an interesting
observable. The random state model \eqref{e:random2} was used in
\cite{BlumGnutSmil02} to predict the following properties:
$$
\frac{\la \tilde\nu(k) \ra}{k}\xrightarrow{k\to\infty}
\frac{|\partial\Omega|}{2\pi}\,,\qquad
\frac{\Var(\tilde\nu(k))}{k}\xrightarrow{k\to\infty} 0.0769\,|\partial\Omega|\,.
$$
The above expectation value was rigorously proved by Toth and
Wigman \cite{TothWig08} within the following random model:
they considered the eigenstates of the Laplacian
$\Delta_\Omega$ on some planar domain $\Omega$, and defined random
linear combinations
in frequency intervals $k_j\in[K,K+1]$:
\be\label{e:random3}
\psi_{rand,K}=\sum_{k_j\in [K,K+1]} a_j\,\psi_j, 
\ee
with the $a_j$ i.i.d. normal Gaussians. 

\subsection{Macroscopic distribution of the nodal set\label{s:nodal-macro}}

The volume (or length) of the nodal set of eigenfunctions is another
interesting quantity. A priori, it should be less sensitive to
perturbations than the nodal count $\nu_j$. Several rigorous results
have been obtained for this quantity.
Donnelly and Fefferman \cite{DonnFeff88} showed that, for any $d$-dimensional
compact real-analytic manifold, the
$(d-1)$-dimensional volume of the nodal set of any Laplacian
eigenstate $\psi_j$ satisfies the bounds
$$
C^{-1}\,k_j \leq \Vol_{d-1} \cN(\psi_j) \leq C\,k_j\,,
$$
for some $C>0$ depending on the manifold. The statistics of this
volume has been investigated for various ensembles of random
states \cite{Ber84,Berry02,RudWig08}. The average volume grows like $c_Mk$, with a
constant $c_M>0$ depending on the manifold. Estimates
for the variance are more difficult to obtain. Berry \cite{Berry02}
argued that for the 2-dimensional random model \eqref{e:random}, the variance
should be of order $\log(k)$, showing an unusually strong concentration property
for this random variable. Such a logarithmic variance was recently proved by
Wigman in the case of random spherical harmonics of the 2-sphere \cite{Wig09}.


\medskip

Counting or volume estimates do not provide any information on
the spatial localization of the nodal set. At the microscopic level,
Br\"uning \cite{Brue78} showed that for any compact riemannian manifold, the nodal
set $\cN(\psi_j)$ is ``dense'' at the scale of the
wavelength: for some constant $C>0$, any ball $B(x,C/k_j)$ intersects $\cN(\psi_j)$.

One can also study the ``macroscopic distribution''
of the zero set, by considering the
$(d-1)$-dimensional riemannian measure on $\cN(\psi_j)$:
\be
\forall f\in C^0(M),\qquad \tilde\mu^Z_{\psi_j}(f)\defeq \int_{\cN(\psi_j)}f(x)\,d\Vol_{d-1}(x)\,.
\ee
Similarly with the case of the density $|\psi_j(x)|^2$ or its phase
space cousins, the spatial
distribution of the nodal set can then be described by the weak-$*$ limits of the
renormalized measures $\mu^Z_{\psi_j}=\frac{\tilde\mu^Z_{\psi_j}}{\tilde\mu^Z_{\psi_j}(M)}$ in the high-frequency limit. 
In the case of chaotic
eigenstates, the following conjecture\footnote{To my knowledge, this conjecture was
  first stated by S.~
Zelditch} seems a reasonable ``dual'' to the QUE conjecture:
\begin{Conj}Let $(M,g)$ be a compact smooth riemannian manifold, with an
  ergodic geodesic flow. Then, for any orthonormal basis
  $(\psi_j)_{j\geq 1}$, the probability measures $\mu^z_{\psi_j}$
  weak-$*$ converge to the Lebesgue measure on $M$, in the limit $j\to\infty$.
\end{Conj}
This conjecture is completely open. One slight
weakening would be to request that the convergence holds for a density
1 subsequence of eigenstates, like in Thm.~\ref{thm:QE}. A similar property can
be proved in a complex analytic setting (see
\S\ref{s:complex-zeros}). For $M$ a real analytic manifold,
eigenfunctions $\psi_j$ can be analytically continued into holomorphic
functions $\psi_j^{\IC}$ in some complex
neighbourhood $M^{\IC}$ of $M$. For 
$(\psi_j)_{j\in S}$ a sequence of {\em ergodic} eigenfunctions (i.e., admitting the Liouville measure as sole semiclassical measure), Zelditch has
shown \cite{Zel07} that the measures associated with the (complex) nodal set of
$\psi_j^{\IC}$ converge to some limit measure on $M^{\IC}$; however, his result says nothing about
the distribution of the real
zeros (that is, of the set $\cN(\psi_j^{\IC})\cap M$).

Once more, it is easier to deal with some class of random states 
than with eigenstates. In the case of random spherical harmonics on the
sphere, the random ensemble is rotation invariant, so for each level
$n$ the expectation of the measure $\mu^Z_{\psi_{rand,n}}$ is equal to the
(normalized) Lebesgue measure.
Zelditch \cite{Zel09} generalized this result to arbitrary compact riemannian
manifolds $M$, by considering random superpositions of
eigenstates of the type \eqref{e:random3}: he showed that the expectation of
$\mu^Z_{\psi_{rand,K}}$ converges to the Legesgue measure in the
limit $K\to\infty$. 

At the moment the study of nodal sets of random wavefunctions represents
lively field of research in probability theory \cite{NazSod10,Wig11},
so the above results are certainly not exhaustive.

\subsection{Nodal sets for eigenstates of quantum maps on the torus}

So far we have only considered the nodal set for the eigenfunctions
$\psi_j$ of the Laplacian, viewed in their spatial representation
$\psi\in L^2(\Omega)$. These eigenfunctions are real valued,
due to the fact that $\Delta$ is a real operator. At the classical level, this reality
corresponds to the fact that the classical flow is time reversal invariant. 

In the case of quantized maps on the torus with time reversal symmetry, the eigenvectors $\psi_j$
also have real position components $\psi_j(\ell/N)=\la
e_\ell,\psi_j\ra$, $\ell=0,\ldots,N-1$. Keating, Mezzadri and Monastra
\cite{KeaMezMon03} have defined the nodal
domains of such states as the intervals
$\{\ell_1,\ell_1+1,\ldots,\ell_2-1\}$ on which $\psi_j(\ell/N)$ has a
constant sign. Using this definition, they computed the exact nodal
statistics for a random state \eqref{e:random-t2} (with real
Gaussian coefficients $a_j$), and numerically showed that these
statistics are satisfied by eigenstates of a generic quantized Anosov
map (namely, a perturbation of the symplectomorphism $S_{DEGI}$).

In a further publication \cite{KeaMarWill06}, these ideas are extended
to chaotic maps on the $4$-dimensional torus. In this setting,
they show that the random state Ansatz directly leads to a 
percolation model on a triangular lattice, with the same
thermodynamic behaviour as the model of \cite{BogoSchm02}. From this
remark, they conjecture that, in the appropriate scaling
limit, the boudaries of the nodal domains for chaotic eigenstates
belong to the universality class of ${\rm SLE}_6$ curves, and support
this claim by some numerics.

\subsection{Husimi nodal sets \label{s:complex-zeros}}
\begin{figure}
\begin{center}
\includegraphics[angle=-90,width=1.\textwidth]{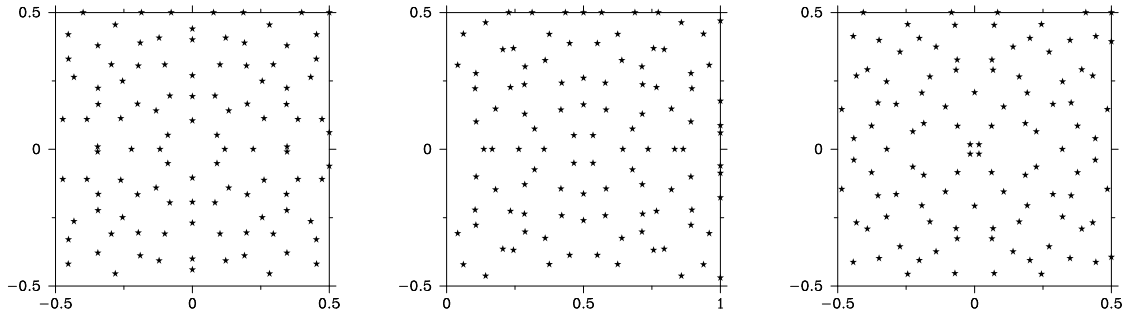}
\caption{Stellar representation for 3 eigenstates of the quantum cat
  map $U_N(S_{DEGI})$, the Husimi densities of which were
  shown in Fig.~\ref{f:cat-Hecke}, top, in a {\em different
    order}. Can you guess
  the permutation? \label{f:cat-3states-stellar}}
\end{center}
\end{figure}
The position representation $\psi(x)$ is physically
natural, but it is not always the most appropriate to investigate semiclassical
properties: phase space representations of the quantum states offer
valuable informations, and are more easily compared with invariant
sets of the classical flow. We have singled out two such phase space
representations: the Wigner function and the
Bargmann-Husimi representation. Both can be defined on $T^*\IR^d$, but also
on the tori $\IT^{2d}$. The Wigner function is real and changes sign,
so its nodal set is an interesting observable. 
However, in this section we will focus on the nodal sets of the
Husimi (or Bargmann) functions.

The Gaussian wavepackets \eqref{e:wavepacket} can be appropriately renormalized into
states
$\tilde\varphi_{\vx}\propto\varphi_{\vx}$ depending {\em antiholomorphically} on the complex
variable $z=x-ip$. As a result, the Bargmann function associated with
$\psi\in L^2(\IR^d)$,
$$
z\mapsto \cB\psi(z)=\la\varphi_{z},\psi\ra\,,
$$
is an entire function of $z$. In dimension $d=1$, the nodal set of
$\cB\psi(z)$ is thus a discrete set of points in $\IC$, which we will
denote by
$\cZ_\psi$. More interestingly,
through Hadamard's factorization one can essentially recover
the Bargmann function $\cB\psi$ from its nodal set, 
and therefore the quantum state $\psi$. Assuming $0\not\in\cZ_\psi$,
we have
$$
\cB\psi(z)=e^{\alpha z^2+\beta z+\gamma}\prod_{0 \neq z_i\in\cZ_\psi}(1-z/z_i)\,e^{\frac{z}{z_i}+\frac12\frac{z^2}{z_i^2}}\,,
$$
leaving only $3$ undetermined parameters.
Leboeuf and Voros \cite{LebVor90} called the set $\cZ_\psi$ the {\em stellar representation}
of $\psi$, and proposed to characterize the chaotic eigenstates
using this representation. 
This idea is especially
appealing in the case of a compact phase space like the 2-torus: in
that case the Bargmann function $\cB\psi$ of a state $\psi_N\in\hn$ is
an entire function on $\IC$ satisfying quasiperiodicity
conditions, so that its nodal set is $\IZ^2$-periodic and contains
exactly $N$ zeros in each fundamental cell. One
can then reconstruct the state $\psi_N$ from this set of $N$ points on
$\t2$ (which we denote by $\cZ^{\t2}_{\psi_N}$):
$$
\cB\psi_N(z)=e^{\gamma}\,\prod_{z_i\in\cZ^{\t2}_{\psi_N}} \theta(z-z_i)\,,
$$
where $\theta(z)$ is a fixed Jacobi theta function vanishing on
$\IZ+i\IZ$, and $\gamma$ a normalization constant. 
This stellar representation is
{\em exact}, {\em minimal} ($N$ complex points represent $\psi\in\hn\equiv
\IC^N$) and lives in {\em phase space}. The conjugation of these three
properties makes it interesting from a semiclassical point of view. 

In \cite{LebVor90} the authors noticed a stark difference
between the nodal patterns of integrable vs. chaotic eigenstates. In
the integrable case, zeros are regularly {\em aligned along certain curves},
which were identified as anti-Stokes lines in the complex WKB
formalism. Namely, the Bargmann function can be approximated by a WKB
Ansatz similar with \eqref{e:Lagrangian} with holomorphic phase functions
$S_j(z)$, and anti-Stokes lines are defined by the equations $\Im
(S_j(z)-S_k(z))=0$ in regions where $e^{iS_j(z)/\hbar}$ and
$e^{iS_k(z)/\hbar}$ dominate the other terms. These anti-Stokes lines
are sitting at the ``antipodes'' of the lagrangian curve on which
the Husimi density is concentrated (for an example of such a state, see Fig.~\ref{f:integ}, center).

On the opposite, the zeros of chaotic eigenstates appear to be
equidistributed across the whole torus (see
Fig.~\ref{f:cat-3states-stellar}), like the Husimi density itself.
This fact was checked on other systems, e.g. planar
billiards, for which the stellar representations of the boundary
functions $\partial_n\psi_j(x)$ were investigated in \cite{TuaVor95}, leading to
similar conclusions. This observation was followed by a rigorous
statement, which we express by defining (using the same notation as
in the previous section) the ``stellar measure'' of a state $\psi_N\in \hn$:
$$
\mu^Z_{\psi_N} \defeq N^{-1}\sum_{z_i\in\cZ^{\t2}_{\psi_N}}\delta_{z_i}\,.
$$
\begin{Theorem}\cite{NonVor98}\label{t:zeros-equid}
Assume that a sequence of normalized states $(\psi_N\in\hn)_{N\geq 1}$
becomes equidistributed on $\t2$ in the limit $N\to\infty$ (that is, their
Husimi measures $\mu_{\psi_N}^H$ weak-$*$ converge to the Liouville measure
$\mu_L$). 

Then, the corresponding stellar measures $\mu^Z_{\psi_N}$ also weak-$*$ converge to $\mu_L$.
\end{Theorem}
Using the quantum ergodicity theorem (or quantum unique ergodicity
when available), one deduces that (almost) all sequences of chaotic
eigenstates have asymptotically equidistributed Husimi nodal sets.

This result was proved independently, and in greater generality, by
Shiffman and Zelditch \cite{ShifZel99}. It was also extended by Rudnick
to the case of holomorphic cusp forms on the modular surface
\cite{Rud05}. The strategy is to first
show that the {\em electrostatic potential}
$$
u_{\psi_N}(\vx) = N^{-1}\log\cH_{\psi_N}(\bx) = 2N^{-1}\log|\cB\psi_N(z)|-\pi|z|^2
$$
decays (in $L^1$) in the semiclassical limit, and then use the fact
that $\mu_{\psi_N}^Z=4\pi\Delta u_{\psi_N}$. This use of potential theory (specific to
the holomorphic setting) explains why such a corresponding statement has not been
proved yet for the nodal set of real eigenfunctions (see the
discussion at the end of \S\ref{s:nodal-macro}).
To my knowledge, this result is the only rigorous one concerning the
stellar representation of chaotic eigenstates. 

In parallel, many studies have been
devoted to the statistical properties of stellar representation of random states
\eqref{e:random-t2}, which can then be compared with those of chaotic
eigenstates. 
Zeros of random entire functions (e.g. random polynomials) have a
long history in probability theory, see e.g. the recent review
\cite{NazSod10}, which mentions the works of Kac,
Littlewood, Offord, Rice. The topic has been revived in the years 1990
through questions appearing in quantum chaos \cite{LebShuk96,BoBoLe96,Han96}. 

For a Gaussian ensemble of random states like \eqref{e:random-t2}, one can explicitly
compute the $n$-point correlation functions of the
zeros: besides being equidistributed, the zeros statistically {\em
  repel} each other quadratically on the
microscopic scale $N^{-1/2}$ (the typical distance between nearby
zeros), but are uncorrelated at larger distances. Such a local
repulsion (which imposes a certain {\em rigidity} of random nodal
sets) holds in great generality, showing a form of
universality at the microscopic scale \cite{BleShiZel00,NazSod10}.

The study of \cite{NonVor98} suggested that the localization
properties of the Husimi measure (e.g. a scar on a periodic orbit) could
not be directly visualized in the distribution of the few zeros near
the scarring orbit, but rather in the {\em collective} distribution of all
zeros. We thus studied in detail the Fourier coefficients of the stellar
measures, 
$$
\mu^Z_{\psi}(e^{2i\pi \vk\cdot\vx}),\qquad 0\neq \vk\in\IZ^2\,.
$$ 
In the case of the random model \eqref{e:random-t2},
the variance of the Fourier coefficients could be explicitly computed:
for fixed
$\vk\neq 0$ and $N\gg 1$, the variance is
$\sim \pi^2\zeta(3)|\vk|^4/N^3$, showing that typical Fourier coefficients are
of size $\sim N^{-3/2}$.
In the case of chaotic eigenstates, we conjectured for each fixed $\vk\neq
0$ an absolute upper bound 
$$\mu^Z_{\psi_{j,N}}(e^{2i\pi \vk\cdot\vx})=o(N^{-1})\,,
$$
while the equidistribution of Thm.~\ref{t:zeros-equid} only forces these
coefficients to be $o(1)$. 
We also argued that the presence of scars in individual eigenfunctions
could be detected through the abnomally large values of a few low-$\vk$ coefficients.

\medskip

To finish this section, let us mention a
few recent rigorous results concerning zeros of random holomorphic
functions, after \cite{NazSod10}:
\begin{itemize}
\item central limit theorems and large
deviation estimates for the ``linear statistics'' $\mu^Z(f)$, with $f$ a fixed test
function (including the characteristic function on a bounded domain)
\item comparison with other {\it point processes}: the Ginibre
  ensemble obtained by the spectrum of random Gaussian complex matrices;
  ensembles of randomly deformed lattices. Here comes the question of the ``best matching'' between a
  random set of points with a square lattice.
\item how to use a set of random zeros to partition the plane into cells.
\end{itemize}
To my knowledge these properties have not been studied in the
framework of Husimi functions of chaotic eigenstates. 

To summarize, the stellar representation provides
complementary ({\em dual}) information to the
macroscopic features of the Husimi measures. The very nonlinear
relation with the wavefunction makes its study difficult, but at the
same time interesting. 


\end{document}